\documentclass[review]{elsarticle}

\usepackage{hyperref}
\usepackage{mathtools}
\usepackage{epsfig}
\usepackage{amsmath}
\usepackage{amsfonts}
\usepackage{amsthm}
\usepackage{amsbsy}
\usepackage{appendix}
\usepackage{epstopdf}
\usepackage[table]{xcolor}
\usepackage{subfloat}
\usepackage{float}
\usepackage[thinlines]{easytable}
\usepackage{array}
\usepackage{color}
\usepackage{inputenc}
\usepackage{hyperref}
\usepackage[ruled,vlined]{algorithm2e}
\usepackage{verbatim}
\usepackage{caption}

\numberwithin{equation}{section}
\allowdisplaybreaks

\theoremstyle{plain}

\newtheorem{assump}{\bf Assumption}[section]

\newtheorem{definition}{\bf Definition}[section]

\theoremstyle{remark}
\newtheorem{remark}{\bf Remark}[section]
\theoremstyle{remark}

\journal{Journal of \LaTeX\ Templates}

\begin{document}

\begin{frontmatter}

\title{Model Reduction 
via Dynamic Mode Decomposition}
\author[mymainaddress]{Hannah Lu\fnref{email1}}
\fntext[email1]{email: \texttt{hannahlu{@}stanford.edu}}
\author[mymainaddress]{Daniel M. Tartakovsky\corref{mycorrespondingauthor}}
\cortext[mycorrespondingauthor]{Corresponding author}
\ead{tartakovsky@stanford.edu}
\address[mymainaddress]{Department of Energy Resources Engineering, Stanford University, Stanford, CA 94305, USA}

\begin{abstract}
This work proposes a new framework of model reduction for parametric complex systems. The framework employs a popular model reduction technique dynamic mode decomposition (DMD), which is capable of combining data-driven learning and physics ingredients based on the Koopman operator theory. In the offline step of the proposed framework, DMD constructs a low-rank linear surrogate model for the high dimensional quantities of interest (QoIs) derived from the (nonlinear) complex high fidelity models (HFMs) of unknown forms. Then in the online step, the resulting local reduced order bases (ROBs) and parametric reduced order models (PROMs) at the training parameter sample points are interpolated to construct a new PROM with the corresponding ROB for a new set of target/test parameter values. The interpolations need to be done on the appropriate manifolds within consistent sets of generalized coordinates. The proposed framework is illustrated by numerical examples for both linear and nonlinear problems. In particular, its advantages in computational costs and accuracy are demonstrated by the comparisons with projection-based proper orthogonal decomposition (POD)-PROM and Kriging.
\end{abstract}


\end{frontmatter}


\section{Introduction}

Physics-based modeling and simulation is playing a significant role across many applications in engineering and science. Despite continuing advances in software and hardware development, including high-performance computing, high fidelity simulations for complex systems remain a demanding and elusive task. That is especially so in sensitivity analysis, uncertainty quantification, inverse modeling, where many simulation runs are required. Real-time applications (e.g., decision making and optimization) also require the parametric computational models, which represent the underlying dynamics in various scenarios, to be simulated accurately within the limited computational time or capacities.

Model reduction techniques can significantly reduce the (prohibitively) high computational costs of physics-based high fidelity simulations, while capturing the essential features of the underlying dynamics. Such techniques have been used extensively in computational engineering and can be grouped in two general classes. The first are physics-based/projection-based approaches (see~\cite{benner2015survey} for a comprehensive review), which seek for a Galerkin projection from a HFM to a representative surrogate model on a  low-dimensional ``trial" subspace for the system states. A prime example of this class is POD~\cite{lumley1967structure,kerschen2005method} and discrete empirical interpolation (DEIM)~\cite{chaturantabut2010nonlinear} for linear and nonlinear systems respectively, both of which are grounded in singular value decomposition (SVD). They obtain the linear trial subspaces determined by the SVD analysis in the time domain. Related approaches in the frequency domain include balanced truncation~\cite{rowley2005model,moore1981principal} and rational interpolation~\cite{gugercin2008h_2,antoulas2020interpolatory}. Alternatively, nonlinear trial manifolds can be obtained using deep convolutional autoencoders in a recent work~\cite{lee2020model}. The computational saving stems from replacing the high-dimensional full  system with its lower-dimensional counterpart for future prediction. Such reduced order models (ROMs) are physics-based in the sense that they inherit the (presumably known) dynamic operator from the projection.

As HFMs become increasingly complex and data becomes more available, there is a growing need for data-driven ROMs, which belong to the second class. Instead of making a direct dimension-reduction of the underlying high-fidelity code that produces the data, data-driven ROMs aim to reduce its complexity by learning the dynamics of the state variables or QoIs directly from the full model’s output and/or observational data. Machine learning techniques, such as Gaussian process regression (also known as Kriging or response surface estimation)~\cite{rasmussen2003gaussian,pau2013reduced}, random forest (RF)~\cite{booker2014comparing,naghibi2016gis}, dynamic mode decomposition (DMD)~\cite{schmid2010dynamic,kutz2016dynamic}, operator inference~\cite{peherstorfer2016data,mcquarrie2021non} and neural networks (NN)~\cite{hesthaven2018non,qin2019data}, are widely used in the constructions of these data-driven ROMs. Physics ingredients can be embedded in the data-driven learning so that both flexibility and robustness can be enhanced even with limited data accessibility~\cite{nathan2018applied,lu2020lagrangian,lu2020prediction,lu2021dynamic,lu2021extended, raissi2019physics,karniadakis2021physics,zhu2019physics,qian2020lift}.

Many complex systems are designed and analyzed by their dependency on the parameters, which account for variations in shape, material, loading, and boundary and initial conditions. Consequently, the ROMs for these parametric complex systems (known as PROMs) are required to be robust with respect to the variations in parameters. For the first class of physics-based/projection-based PROMs, this issue was discussed in~\cite{epureanu2003parametric,homescu2005error,serban2007effect,lieu2007adaptation} and later resolved in~\cite{amsallem2008interpolation,amsallem2011online} using appropriate interpolations of the ROBs and PROMs on suitable manifolds respectively. Later improvements in interpolation include~\cite{son2013real,zimmermann2014locally,zhang2021gaussian}. The major challenge of the first class PROMs for complex systems is that they are ``physics hungry"  -- the projections can only be designed for the full state variables with full knowledge of the dynamic operators. It is inefficient and sometimes unfeasible to deal with a large number of state variables with their extremely high dimension after spatial discretization and huge systems of (nonlinear) dynamic operators. On the other hand, the second class data-driven PROMs can be ``data hungry". The time evolution of the QoIs, expressed as a function of the parameters and simulated predictors, can present high nonlinearity or/and bifurcations. It becomes challenging for conventional approximators like Gaussian Regression Process to capture the dynamics precisely as discussed in~\cite{amsallem2011online} and verified in our numerical example section~\ref{sec:NS}. Modern nonlinear machine learning methods like RF and NN show advantageous performances in this case but only when enough amount of high-fidelity data are available for training~\cite{lu2021data,bhattacharya2021model, sentz2021reduced}. Alternative strategies include operator inference~\cite{mcquarrie2021non}, which needs assumptions of polynomial structures on the underlying dynamics. 

Given the aforementioned challenges, adapting the data-driven model reduction strategies to the parametric setting is an active area of research. In this work, we propose a physics-aware data-driven DMD-PROM framework, which combines the advantages from both classes of PROMs. The proposed framework consists of two steps: 1). The offline step constructs a low-rank DMD-based (linear) ROM for the observables of the QoIs at each training sample parameter point. The physics-guided selection of observables, which is based on the Koopman operator theory, provides not only better accuracy in the local ROMs, but also a bridge between the understanding of data and physics. 2). The online step constructs a new PROM with corresponding ROB for each target/test parameter point by interpolating the training PROMs and ROBs. The interpolations are done on suitable manifolds in the same fashion as in the projection-based PROMs. Our framework alleviates the constraints on the access to the HFMs in projection-based PROMs and the need for large amount of data in pure data-driven PROMs. 

The remainder of the paper is organized as follows. The problem of interest is formulated in section~\ref{sec:2}. Then we establish the general methodology and summarize the practical algorithms in section~\ref{sec:3}. In section~\ref{sec:4}, several numerical examples are presented to demonstrate the accuracy and robustness of the proposed framework, with comparison with other methods. Key results, their implication for applications, and challenges and future work are summarized in section~\ref{sec:5}.

\section{Problem Formulation}\label{sec:2}

We consider a complex system that is described by $N_\text{sv}'$ state variables $\mathbf s(\mathbf x,t;\mathbf p) = \{s_1,\dots,s_{N_\text{sv}'} \}$, varying in space $\mathbf x \in \mathcal D$, time $t \in [0,T]$ and parameter $\mathbf p\in \Omega$ throughout the simulation domain $\mathcal D$ during the simulation time interval $[0,T]$ within the parameter range $\Omega$. The spatiotemporal evolution of these state variables is described by a system of coupled partial-differential equations (PDEs)
\begin{align}\label{eq:genPDF}
    \frac{\partial s_i}{\partial t} = \phi_i(\mathbf s; \mathbf p), \quad (\mathbf x,t) \in \mathcal D \times (0,T]; \qquad i = 1,\dots,N_\text{sv}',
\end{align}
where $\phi_i$ are (linear/nonlinear) differential operators that contain spatial derivatives, and $\mathbf p = \{p_1,\dots,$ $p_{N_\text{par}} \}$ is a set of $N_\text{par}$ parameters. Problems of this kind have to be solved numerically, which requires a discretization of the spatial domain $\mathcal D$ into $N_\text{el}$ elements/nodes, leading to a discretized state variable $\mathbf S(t;\mathbf p)$ of very high dimension $N_\text{sv} = N_\text{sv}'\times N_\text{el}$, satisfying
\begin{equation}\label{eq:ODE}
\frac{d \mathbf S}{dt}(t;\mathbf p) = \boldsymbol \Phi (\mathbf S(t;\mathbf p);\mathbf p),
\end{equation}
where $\boldsymbol \Phi$ are linear/nonlinear functions.

More often than not, this model output has to be post-processed to compute $N_Q$ QoIs $\mathbf Q(t;\mathbf p) = [Q_1(t;\mathbf p),\dots,Q_{N_Q}(t;\mathbf p) ]^\top$, such that 
\begin{equation}\label{eq:QoI}
Q_i (t;\mathbf p)= \mathcal M_i(\mathbf S(t;\mathbf p)),\quad \text{for}\quad i=1,\dots,N_Q.
\end{equation}
 The maps $\mathcal M_i$ can represent, e.g., numerical approximations of the integrals over $\mathcal D$ or a streamline. The QoIs are usually much easier to visualize and comprehend than the raw output $\mathbf S (t;\mathbf p)$. We are interested in the HFM for the QoIs
 \begin{equation}\label{eq:HFM}
 \frac{d\mathbf Q}{dt}(t;\mathbf p) = \mathbf F(\mathbf Q(t;\mathbf p);\mathbf p)
 \end{equation}
instead of  the complex system~\eqref{eq:ODE} for $\mathbf S(t;\mathbf p)$. Notice that the function $\mathbf F$ has a complicated form which depends on the knowledge of the original system~\eqref{eq:ODE}. The corresponding discrete-time dynamical system is described by
\begin{equation}\label{eq:HFM-dis}
\mathbf Q(t_{k+1};\mathbf p) = \mathbf F_{\Delta t}(\mathbf Q(t_{k};\mathbf p);\mathbf p) := \mathbf Q(t_{k};\mathbf p)+\int_{t_k}^{t_k+\Delta t}\mathbf F(\mathbf Q(\tau;\mathbf p))d\tau
\end{equation}
for the uniform time discretization $t_k = k\Delta t\in [0,T], k = 0,\cdots, N_T$.


Given $N_\text{snap}$ snapshots of the trajectories of the QoIs computed from~\eqref{eq:ODE} and~\eqref{eq:QoI} at sampled parameter points $\{\mathbf p^{(1)},\cdots, \mathbf p^{(N_\text{MC})}\}$,
\begin{equation}\label{eq:data}
\{\mathbf Q(t_0;\mathbf p^{(1)}),\cdots,\mathbf Q(t_{N_\text{snap}};\mathbf p^{(1)}),\cdots,\mathbf Q(t_0;\mathbf p^{(N_\text{MC})}),\cdots,\mathbf Q(t_{N_\text{snap}};\mathbf p^{(N_\text{MC})})\}, \quad N_\text{snap}\leq N_T,
\end{equation}
we aim to construct a surrogate model of reduced dimension for the discretized HFM~\eqref{eq:HFM-dis}. This will allow us to predict the trajectory of the QoIs at an unsampled parameter point $\mathbf p^*\notin \{\mathbf p^{(1)},\cdots, \mathbf p^{(N_\text{MC})}\}$, i.e., $\{\mathbf Q(t_0;\mathbf p^*),\cdots,\mathbf Q(t_{N_T};\mathbf p^*)\}$, directly at a low cost without computing the complex high dimensional system~\eqref{eq:ODE} and~\eqref{eq:QoI}.

\section{Methodology}\label{sec:3}
We propose a DMD-PROM framework consisting of the following offline step and online step in section~\ref{sec:offline} and section~\ref{sec:online} respectively.

\subsection{Offline Step: DMD-based Surrogate Models}\label{sec:offline}
Without the accessibility to the original system~\eqref{eq:ODE} and~\eqref{eq:QoI}, one can regard the dynamic function $\mathbf F$ as unknown in~\eqref{eq:HFM} (and $\mathbf F_{\Delta t}$ unknown in~\eqref{eq:HFM-dis} correspondingly). The goal of this step is to construct a surrogate linear model of reduced dimension for the unknown discrete HFM~\eqref{eq:HFM-dis} from the dataset~\eqref{eq:data}. 

We first review the Koopman operator theory, which allows one to handle the potential nonlinearity in the unknown dynamic $\mathbf F$ and $\mathbf F_{\Delta t}$:
\begin{definition}[Koopman operator~\cite{kutz2016dynamic}]
For nonlinear dynamic system~\eqref{eq:HFM}, the Koopman operator $\mathcal K^{\mathbf p}$ is an infinite-dimensional linear operator that acts on all observable functions $g: \mathcal M\to \mathbb C$ so that
\begin{equation}\label{eq:Koopman}
\mathcal K^{\mathbf p} g(\mathbf Q(t;\mathbf p)) = g(\mathbf F(\mathbf Q(t;\mathbf p);\mathbf p)).
\end{equation}
For discrete dynamic system~\eqref{eq:HFM-dis}, the discrete-time Koopman operator $\mathcal K_{\Delta t}^{\mathbf p}$ is
\begin{equation}\label{eq:Koopman-dis}
\mathcal K_{\Delta t} ^{\mathbf p}g(\mathbf Q(t_{k+1};\mathbf p)) = g(\mathbf F_{\Delta t}(\mathbf Q(t_k;\mathbf p);\mathbf p)) = g(\mathbf Q(t_{k+1};\mathbf p)).
\end{equation}
\end{definition}

The Koopman operator transforms the finite-dimensional nonlinear problem~\eqref{eq:HFM} in the state space into the infinite-dimensional linear problem~\eqref{eq:Koopman} in the observable space. Since $\mathcal K_{\Delta t}^{\mathbf p}$ is an infinite-dimensional linear operator, it is equipped with infinite eigenvalues $\{\lambda_k(\mathbf p)\}_{k=1}^\infty$ and eigenfunctions $\{\phi_k(\mathbf p)\}_{k=1}^\infty$. In practice, one has to make a finite approximation of the eigenvalues and eigenfunctions. The following assumption is essential to both a finite-dimensional approximation and the choice of observables.

\begin{assump}\label{assump1}
Let $\mathbf y (t_k;\mathbf p)$ denote a $N\times 1$ vector of observables,
\begin{equation}
\mathbf y(t_k;\mathbf p) = \mathbf g(\mathbf Q(t_k;\mathbf p) )= \begin{bmatrix}
g_1(\mathbf Q(t_k;\mathbf p))\\
\vdots\\
g_N(\mathbf Q(t_k;\mathbf p))
\end{bmatrix},
\end{equation}
where $g_j:\mathcal M\to \mathbb C$ is an observable function with $j = 1,...,N$. If the chosen observable $\mathbf g$ is restricted to an invariant subspace spanned by eigenfunctions of the Koopman operator $\mathcal K_{\Delta t}^{\mathbf p}$, then it induces a linear operator $\mathbf K(\mathbf p)$ that is finite-dimensional and advances these eigenobservable functions on this subspace~\cite{brunton2016koopman}.
\end{assump}
Based on Assumption~\ref{assump1}, one can deploy the DMD algorithm~\ref{alg:DMD}. to approximate the $N$-dimensional linear operator $\mathbf K(\mathbf p)$ and its low dimensional approximation $\mathbf K_r(\mathbf p)$ of rank $r$. At each parameter point $\mathbf p^{(i)}, i = 1,\cdots, N_{\text{MC}}$, the discrete HFM~\eqref{eq:HFM-dis} on state space is approximated by a $N$-dimensional linear surrogate model
\begin{equation}\label{eq:surrogate}
\mathbf y(t_{k+1};\mathbf p^{(i)}) = \mathbf K(\mathbf p^{(i)})\mathbf y(t_{k};\mathbf p^{(i)})
\end{equation}
on observable space. The two spaces are connected  by the observable function $\mathbf g$ and its inverse $\mathbf g^{-1}$. Algorithm~\ref{alg:DMD} induces the following ROM for~\eqref{eq:surrogate} directly,
\begin{equation}
\label{eq:ROM}
\mathbf q(t_{k+1};\mathbf p^{(i)}) = \mathbf K_r(\mathbf p^{(i)})\mathbf q(t_{k};\mathbf p^{(i)}),
\end{equation}
where 
\begin{equation}
\mathbf y(t_{k};\mathbf p^{(i)}) = \mathbf V (\mathbf p^{(i)})\mathbf q(t_{k};\mathbf p^{(i)}),\quad  \mathbf K_r(\mathbf p^{(i)}) =  \mathbf V (\mathbf p^{(i)}) ^\top \mathbf K(\mathbf p^{(i)})\mathbf V (\mathbf p^{(i)}).
\end{equation}

\begin{algorithm}[H]
Input: $\{\mathbf Q(t_0;\mathbf p^{(i)}),\cdots,\mathbf Q(t_{N_\text{snap}};\mathbf p^{(i)})\}$, observable function $\mathbf g$,
\begin{enumerate}
\item Create the data matrices of observables
\begin{equation}
\mathbf Y_1(\mathbf p^{(i)}) = \begin{bmatrix}
|&&|\\
\mathbf y(t_0;\mathbf p^{(i)})&\cdots&\mathbf y(t_{N_\text{snap}-1};\mathbf p^{(i)})\\
|&&|
\end{bmatrix}, \quad \mathbf Y_2 (\mathbf p^{(i)}) = \begin{bmatrix}
|&&|\\
\mathbf y(t_1;\mathbf p^{(i)})&\cdots&\mathbf y(t_{N_\text{snap}};\mathbf p^{(i)})\\
|&&|
\end{bmatrix},
\end{equation}
where each column is given by $\mathbf y(t_k; \mathbf p^{(i)}) = \mathbf g(\mathbf Q(t_k;\mathbf p^{(i)}))$.
\item Apply SVD $\mathbf Y_1(\mathbf p^{(i)}) \approx \mathbf V(\mathbf p^{(i)}) \boldsymbol \Sigma(\mathbf p^{(i)}) \mathbf Z (\mathbf p^{(i)}) ^*$ with $\mathbf V(\mathbf p^{(i)}) \in \mathbb C^{N\times r}$, $\boldsymbol \Sigma(\mathbf p^{(i)})  \in \mathbb C^{r\times r}$, $\mathbf Z(\mathbf p^{(i)}) \in \mathbb C^{r\times (N_\text{snap}-1)}$, where $r$ is the truncated rank chosen by certain criteria and should be the same for all $i = 1,\cdots, N_\text{MC}$.
\item Compute $\mathbf K_r(\mathbf p^{(i)}) = \mathbf V(\mathbf p^{(i)})^*\mathbf Y_2(\mathbf p^{(i)})\mathbf Z(\mathbf p^{(i)})\boldsymbol \Sigma(\mathbf p^{(i)})^{-1}$ as an $r\times r$ low-rank approximation for $\mathbf K(\mathbf p^{(i)})$.
\item Compute $\mathbf P^{(i,j)} = \mathbf V(\mathbf p^{(i)})^\top\mathbf V(\mathbf p^{(j)})$ for $i,j = 1,\cdots, N_\text{MC}$.
\end{enumerate}
Output: $\mathbf V (\mathbf p^{(i)})$, $\mathbf K_r(\mathbf p^{(i)})$ and $\mathbf P^{(i,j)}$.
\caption{DMD algorithm on observable space~\cite{kutz2016dynamic} for parameter point $\mathbf p^{(i)}, i = 1,\cdots, N_{\text{MC}}$.}
\label{alg:DMD}
\end{algorithm}

\begin{remark}
Notice that the above step can be done offline. $\mathbf P^{(i,j)}$ are precomputed for the later online step. Although the offline step takes a majority of the computational costs in the whole framework mostly due to the computation for the high-fidelity training data~\eqref{eq:data}, the output of this step can be pre-computed and stored efficiently (the output storage is $\left(N\times r+r\times r+r\times r\times \frac{N_\text{MC}+1}{2}\right)\times N_\text{MC}$). In practice, real-time applications can be enabled as long as the online step is sufficient enough.
\end{remark}

\begin{remark}
Connections between the DMD theory and the Koopman spectral analysis under specific conditions on the observables and collected data are established by a theorem in~\cite{tu2013dynamic}. This theorem indicates that judicious selection of the observables is critical to success of the Koopman method. In general, there is no principled way to select observables without expert knowledge of a dynamical system. Machine learning techniques can be deployed to identify relevant terms in the dynamics from data, which guide selection of the observables~\cite{schmidt2009distilling,wang2011predicting,wang2011predicting}. In our numerical examples, we select the observables based on knowledge about the underlying physics (similar to previous works~\cite{lu2020lagrangian,lu2020prediction,lu2021dynamic,lu2021extended,williams2014kernel,li2017extended}).
\end{remark}

\subsection{Online Step: Interpolation of ROBs and PROMs}\label{sec:online}
For an unsampled parameter point $\mathbf p^*\notin \{\mathbf p^{(1)},\cdots, \mathbf p^{(N_\text{MC})}\}$, the goal is to compute $\{\mathbf Q(t_1;\mathbf p^*),\cdots,\mathbf Q(t_{N_T};\mathbf p^*)\}$ via the following PROM without computing the complex high dimensional system~\eqref{eq:ODE} and~\eqref{eq:QoI}:
\begin{equation}\label{eq:ROM-iter}
\mathbf q(t_{k+1};\mathbf p^{*}) = \mathbf K_r(\mathbf p^{*})\mathbf q(t_{k};\mathbf p^{*}).
\end{equation}
Subsequently, the state variable can be computed by
\begin{equation}\label{eq:ROM-state}
\mathbf y(t_{k};\mathbf p^*) = \mathbf V (\mathbf p^*)\mathbf q(t_{k};\mathbf p^*),\quad
\mathbf Q(t_k; \mathbf p^*) = \mathbf g^{-1}(\mathbf y(t_k;\mathbf p^*)).
\end{equation}
 Therefore, the online task includes 1. Computing $\mathbf V (\mathbf p^*)$; 2. Computing $\mathbf K_r (\mathbf p^*)$; 3. Computing $\mathbf q(t_k;\mathbf p^*)$.

\subsubsection{Interpolating ROBs}\label{sec:interpROB}
In order to compute the ROB $\mathbf V (\mathbf p^*)$ needed in~\eqref{eq:ROM-state}, the following interpolation on Grassman manifold is made from $\{ \mathbf V (\mathbf p^{(1)}),\cdots,\mathbf V (\mathbf p^{(N_\text{MC})})\}$. Here we briefly review the interpolation approach proposed in~\cite{amsallem2008interpolation}.

\begin{definition}First we denote the following matrix manifold of interest:
\begin{itemize}
\item Grassmann manifold $\mathcal G(r,N))$ is the set of all subspaces in $\mathbb R^N$  of dimension $r$;
\item Orthogonal Stiefel Manifold $\mathcal S\mathcal T(r,N)$ is the set of orthogonal ROB matrices in $\mathbb R^{r\times N}$. 
\end{itemize}
\end{definition}

Recall that $\mathbf V (\mathbf p^{(i)}) \in \mathbb R^{N\times r}, i = 1,\cdots, N_\text{MC}$, where $r\leq N$, denote the full-rank column matrix, whose columns provide a basis of subspace $\mathcal S_i$ of dimension $r$ in $\mathbb R^N$. The associated ROM is not uniquely defined by the ROB but uniquely defined by the subspace $\mathcal S_i$. Therefore, the correct entity to interpolate is the subspaces $\mathcal S_i$ instead of the ROB  $\mathbf V (\mathbf p^{(i)})$. The goal is now shifted to compute $\mathcal S_* = \text{range}(\mathbf V(\mathbf p^*))$ by interpolating $\{\mathcal S_i\}_{i=1}^{N_{\text{MC}}}$ with access to a ROB $\mathbf V(\mathbf p^*)$.

The subspaces $\mathcal S_i$ belongs to the Grassmann manifold~\cite{absil2004riemannian,boothby2003introduction,helgason2001differential,rahman2005multiscale,edelman1998geometry} $\mathcal G(r,N)$. Each $r$-dimensional subspace $\tilde {\mathcal S}$ of $\mathbb R^N$ can be regarded as a point of $\mathcal G(r,N)$ and nonuniquely represented by a matrix $\tilde{\mathbf V}\in \mathbb R^{N\times r}$ whose columns span the subspace $\tilde{\mathcal S}$. The matrices $\tilde{\mathbf V}$ are chosen among those whose columns form a set of orthonormal vectors of $\mathbb R^N$ and belong to the orthogonal Stiefel manifold $\mathcal S\mathcal T(r,N)$~\cite{absil2004riemannian,edelman1998geometry}. There exists a projection map~\cite{absil2004riemannian} from $\mathcal S\mathcal T(r,N)$ to $\mathcal G(r,N)$, as each matrix in $\mathcal S\mathcal T(r,N)$ spans a subspace of $\mathbb R^N$ of dimension $r$ and different matrices can span the same subspace. At each point $\tilde{\mathcal S}$ of the manifold $\mathcal G(r,N)$, there exists a tangent space~\cite{absil2004riemannian,edelman1998geometry} of the same dimension~\cite{edelman1998geometry}. This space is denoted by $\mathcal T_{\tilde{\mathcal S}}$ and each of its points can be represented by a matrix $\tilde {\boldsymbol \Gamma}\in \mathbb R^{N\times r}$. This tangent space is a vector space where usual interpolation is allowed for the matrices representing points of $\mathcal T_{\tilde{\mathcal S}}$. Let $\boldsymbol \Gamma^i = m_{\mathbf V}(\mathbf V(\mathbf p^{(i)}))$, where $m_{\mathbf V}$ denote  the map from the matrix manifolds $\mathcal G(r,N)$ to the tangent space. Then the goal is to compute $\boldsymbol \Gamma^*$ by performing usual interpolation on $\{\boldsymbol \Gamma^i\}_{i=1}^{N_{\text{MC}}}$ and get back to $\mathbf V(\mathbf p^*)$ through the inverse map $m_{\mathbf V}^{-1}(\boldsymbol \Gamma^*)$.

$m_{\mathbf V}$ is chosen to be the logarithm mapping (which maps the Grassmann manifold to its tangent space), and $m_{\mathbf V}^{-1}$ is chosen to be the exponential mapping (which maps the tangent space to the Grassmann manifold itself). This choice borrows concepts of geodesic path on a Grassmann manifold from differential geometry~\cite{absil2004riemannian,boothby2003introduction,wald1984general,do1992riemannian}. We summarize the algorithm as follows and refer the readers to~\cite{amsallem2008interpolation} for the detailed construction.

\begin{algorithm}[H]
Input: $\{ \mathbf V (\mathbf p^{(i)})\}_{i=1}^{N_{\text{MC}}}$, $\{ \mathbf P^{(i,j)}\}_{i,j=1}^{N_{\text{MC}}}$, $\{\mathbf p^{(i)}\}_{i=1}^{N_{\text{MC}}}$ and target parameter point $\mathbf p^*$,
\begin{enumerate}
\item Denote $\mathcal S_i = \text{range}(\mathbf V(\mathbf p^{(i)}))$, $i = 1,\cdots N_{\text{MC}}$. A point $\mathcal S_{i_0}, i_0\in \{1,\cdots N_{\text{MC}}\}$ of the manifold is chosen as a reference and origin point for interpolation;
\item Select the points $\mathcal S_i, i\in \mathcal I_{i_0}\subset \{1,\cdots, N_{\text{MC}}\}$ which lie in a sufficiently small neighborhood of $\mathcal S_{i_0}$ and map those $\{\mathcal S_i\}_{i\in \mathcal I_{i_0}}$ to matrices $\{\boldsymbol \Gamma^i\}_{i\in \mathcal I_{i_0}}$ representing corresponding points of $\mathcal T_{\mathcal S_{i_0}}$ using the logarithm map $\log_{\mathcal S_{i_0}}$. This can be computed by
\begin{equation}
\begin{aligned}
&(\mathbf I -\mathbf V(\mathbf p^{(i_0)})\mathbf V(\mathbf p^{(i_0)})^\top)\mathbf V(\mathbf p^{(i)})(\mathbf P^{(i_0,i)})^{-1} = \mathbf U_i\boldsymbol \Omega_i\mathbf W_i^\top,\quad \text{(thin SVD)}\\
&\boldsymbol \Gamma^i = \mathbf U_i\tan^{-1}(\boldsymbol \Omega_i)\mathbf W_i^\top.
\end{aligned}
\end{equation}
\item Compute $\boldsymbol \Gamma^*$ by interpolating $\{\boldsymbol \Gamma^i\}_{i\in \mathcal I_{i_0}}$ entry by entry:
\begin{equation}\label{eq:entry-by-entry}
\Gamma_{ij}^* = \mathcal P(\mathbf p^*;\{\Gamma_{ij}^i,\mathbf p^{(i)}\}_{i\in \mathcal I_{i_0}}), \quad i =1,\cdots, N,\quad j = 1,\cdots, r. 
\end{equation}
\item Map the matrix $\boldsymbol \Gamma^*$ representing a point of $\mathcal T_{\mathcal S_{i_0}}$ to the desired subspace $\mathcal S_*$ on $\mathcal G(r,N)$ spanned by a ROB $\mathbf V(\mathbf p^*)$ using the exponential map $\exp_{\mathcal S_{i_0}}$. This can be computed by
\begin{equation}
\begin{aligned}
&\boldsymbol \Gamma^* =  \mathbf U_*\tan^{-1}(\boldsymbol \Omega_*)\mathbf W_*^\top,\quad \text{(thin SVD)}\\
&\mathbf V(\mathbf p^*) = \mathbf V(\mathbf p^{(i_0)})\mathbf W_*\cos(\boldsymbol \Omega_*)+\mathbf U_*\sin(\boldsymbol \Omega_*)
\end{aligned}
\end{equation}
\end{enumerate}
Output: $\mathbf V(\mathbf p^*)$.
\caption{Interpolation of ROBs~\cite{amsallem2008interpolation}}
\label{alg:ROB}
\end{algorithm}

\begin{remark}
The choice of the interpolation method $\mathcal P$ depends on the  dimension of the parameter $N_\text{par}$. When $N_\text{par}=1$, a univariate (typically, a Lagrange type) interpolation method is chosen. Otherwise, a multivariate interpolation scheme (see, for example,~\cite{spath1995one,de1992computational}) is chosen.
\end{remark}

\begin{remark}
Because the logarithmic map $\log_{\mathcal S_{i_0}}$ is defined in a neighborhood of $\mathcal S_{i_0}\in \mathcal G(r,N)$, the method is not sensitive, in principle, to the choice of the reference point $\mathcal S_{i_0}$ in step 1 of Algorithm~\ref{alg:ROB}. This is also confirmed in practice~\cite{amsallem2008interpolation}.
\end{remark}

\subsubsection{Interpolating PROMs}
In order to compute the reduced order operator $\mathbf K_r (\mathbf p^*)$ needed in~\eqref{eq:ROM-iter}, the following interpolation on matrix manifold is made from  ROMs $\{ \mathbf K_r (\mathbf p^{(1)}),\cdots,\mathbf K_r (\mathbf p^{(N_\text{MC})})\}$. We briefly review the approach proposed in~\cite{amsallem2011online} consisting of the following two steps: 
\begin{itemize}
\item Step A). As discussed before, a given ROM can be expressed in a variety of equivalent ROBs. The resulting ROMs may not have been precomputed in the same generalized coordinates system. The validity of an interpolation may crucially depend on the choice of the representative element within each equivalent class. Given the precomputed ROMs $\{ \mathbf K_r (\mathbf p^{(1)}),\cdots,\mathbf K_r (\mathbf p^{(N_\text{MC})})\}$, a set of congruence transformations is determined so that a representative element of the equivalent ROBs for each precomputed ROM is chosen to align the precomputed ROMs into consistent sets of generalized coordinates. The consistency can be enforced by solving the following classical orthogonal Procrustes problems~\cite{van1996matrix}:
\begin{equation}
\min_{\mathbf Q_i,\mathbf Q_i^\top\mathbf Q_i = \mathbf I_r}\|\mathbf V(\mathbf p^{(i)})^\top\mathbf Q_i-\mathbf V(\mathbf p^{(i_0)})\|_F, \forall i = 1,\cdots, N_\text{MC},
\end{equation}
where $i_0\in \{1,\cdots N_{\text{MC}}\}$ is chosen as a reference configuration.
The representative element can be identified by solving the above problem analytically. We summarize the algorithm below and refer the readers to~\cite{amsallem2011online} for details.

\begin{algorithm}[H]
Input: $\{ \mathbf K_r (\mathbf p^{(1)}),\cdots,\mathbf K_r (\mathbf p^{(N_\text{MC})})\}$, $\{ \mathbf P^{(i,j)}\}_{i,j=1}^{N_{\text{MC}}}$, reference configuration choice $i_0$,

\textbf{For}  $i\in \{1,\cdots N_{\text{MC}}\}\setminus\{i_0\}$
\begin{itemize}
\item Compute $\mathbf P^{(i,i_0)} = \mathbf U_{i,i_0}\boldsymbol \Sigma_{i,i_0}\mathbf Z_{i,i_0}^\top$ (SVD),
\item Compute $\mathbf Q_i = \mathbf U_{i,i_0}\mathbf Z_{i,i_0}^\top$,
\item Transform $\tilde{\mathbf K}_r (\mathbf p^{(i)}) = \mathbf Q_i^\top \mathbf K_r (\mathbf p^{(i)}) \mathbf Q_i$
\end{itemize}
\textbf{End}

Output: $\{ \tilde{\mathbf K}_r (\mathbf p^{(1)}),\cdots,\tilde{\mathbf K}_r (\mathbf p^{(N_\text{MC})})\}$.
\caption{Step A of interpolating PROMs~\cite{amsallem2008interpolation}}
\label{alg:stepA}
\end{algorithm}

\begin{remark}
An optimal choice of the reference configuration $i_0$, if it exists, remains an open problem.
\end{remark}
\begin{remark}
The above step A is related to mode tracking procedures based on the modal assurance criterion (MAC)~\cite{ewins2009modal}. The connection is illustrated in~\cite{amsallem2011online}.
\end{remark}

\item Step B). The transformed ROMs $\{ \tilde{\mathbf K}_r (\mathbf p^{(1)}),\cdots,\tilde{\mathbf K}_r (\mathbf p^{(N_\text{MC})})\}$ are interpolated to compute a ROM $\mathbf K_r (\mathbf p^*)$. Similar to section~\ref{sec:interpROB}, this interpolation must be performed on a specific manifold  containing $\{ \tilde{\mathbf K}_r (\mathbf p^{(1)}),\cdots,\tilde{\mathbf K}_r (\mathbf p^{(N_\text{MC})})\}$ and $\mathbf K_r (\mathbf p^*)$ so that the distinctive properties (e.g., orthogonality, nonsingularity) can be preserved. The main idea again is to first map all precomputed matrices to the tangent space to the matrix manifold of interest at a chosen reference point using the logarithm mapping, interpolate the mapped data in this linear vector space, and finally map the interpolated result back to the manifold of interest using the associated exponential map. The algorithm is described as follows:

\begin{algorithm}[H]
Input: $\{ \tilde{\mathbf K}_r (\mathbf p^{(1)}),\cdots,\tilde{\mathbf K}_r (\mathbf p^{(N_\text{MC})})\}$, reference configuration choice $i_0$,
\begin{enumerate}
\item \textbf{For}  $i\in \{1,\cdots N_{\text{MC}}\}\setminus\{i_0\}$
\begin{itemize}
\item Compute $\boldsymbol \Gamma^i = \log_{\tilde{\mathbf K}_r(\mathbf p^{(i_0)})}(\tilde{\mathbf K}_r(\mathbf p^{(i)}))$,
\end{itemize}
\textbf{End}
\item Compute $\boldsymbol \Gamma^*$ by interpolating $\{\boldsymbol \Gamma^i\}_{i\in \mathcal I_{i_0}}$ entry by entry the same as~\eqref{eq:entry-by-entry}.
\item Compute $\mathbf K_r(\mathbf p^*) = \exp_{\tilde{\mathbf K}_r(\mathbf p^{(i_0)})}(\boldsymbol \Gamma^*)$.
\end{enumerate}
Output: $\mathbf K_r(\mathbf p^*)$.
\caption{Step B of interpolating PROMs~\cite{amsallem2008interpolation}}
\label{alg:stepB}
\end{algorithm}
\begin{remark}
The $\log$ and $\exp$ denote the matrix logarithm and exponential respectively. The specific expressions of different matrix manifolds of interest are listed in Table 4.1 of~\cite{amsallem2011online}.
\end{remark}
\end{itemize}
\subsubsection{Iteration-free Computation of the Solution}
With $\mathbf K_r(\mathbf p^*)$ and $\mathbf V(\mathbf p^*)$, one can compute the solution using the iteration-free feature of DMD framework. The algorithm is summarized as below:

\begin{algorithm}[H]
Input: $\mathbf K_r(\mathbf p^*)$, $\mathbf V(\mathbf p^*)$, $\mathbf y(t=0; \mathbf p^*)$, $\mathbf g$.
\begin{enumerate}
\item Compute the eigen-decomposition of $\mathbf K_r(\mathbf p^*)$:
\begin{equation}
\mathbf K_r(\mathbf p^*) \boldsymbol \Psi (\mathbf p^*) = \boldsymbol \Psi (\mathbf p^*) \boldsymbol \Lambda(\mathbf p^*),
\end{equation}    
where columns of $\boldsymbol \Psi (\mathbf p^*)$ are eigenvectors and $ \boldsymbol \Lambda(\mathbf p^*)$ is a diagonal matrix containing the corresponding eigenvalues $\lambda_i, i = 1,\cdots, r$.
\item Compute the DMD modes 
\begin{equation}
\boldsymbol \Phi (\mathbf p^*) = \mathbf V(\mathbf p^*) \boldsymbol \Psi (\mathbf p^*). 
\end{equation}
\item Reconstruct the observables:
\begin{equation}
\mathbf y_\text{DMD}(t_k;\mathbf p^*) = \boldsymbol \Phi (\mathbf p^*)\boldsymbol \Lambda(\mathbf p^*)^k \left(\boldsymbol \Phi (\mathbf p^*)^{-1}\mathbf y(0;\mathbf p^*) \right), k = 1,\cdots, N_T.
\end{equation}
\item Map the observables back to state space:
\begin{equation}
\mathbf Q_\text{DMD}(t_k;\mathbf p^*) = \mathbf g^{-1}(\mathbf y_\text{DMD}(t_k;\mathbf p^*) ).
\end{equation}
\end{enumerate}
Output: $\mathbf Q_\text{DMD}(t_k;\mathbf p^*)$.
\caption{DMD Reconstruction}
\label{alg:dmd-reconst}
\end{algorithm}

\subsection{Algorithm Summary}
The proposed DMD-PROM framework is summarized below:

\begin{algorithm}[H]
\textit{Offline Step:} 

\textbf{For} $i = 1,\cdots,N_\text{MC}$,
$$
\begin{aligned}
&\text{Compute the high fidelity training data~\eqref{eq:data}, } \\
&\text{Input:} \ \{\mathbf Q(t_1;\mathbf p^{(i)}),\cdots,\mathbf Q(t_{N_\text{snap}};\mathbf p^{(i)})\} \ \text{and} \ \mathbf g\xrightarrow{\text{Algorithm~\ref{alg:DMD}}} \text{Output:} \ \mathbf V(\mathbf p^{(i)}),  \mathbf K_r(\mathbf p^{(i)}) \ \text{and} \ \mathbf P^{(i,j)}
\end{aligned}$$

\textbf{End}

\textit{Online Step:}
\begin{itemize}
\item Interpolation of ROBs:
$$\text{Input:} \ \{ \mathbf V (\mathbf p^{(i)})\}_{i=1}^{N_{\text{MC}}}, \ \{ \mathbf P^{(i,j)}\}_{i,j=1}^{N_{\text{MC}}}, \ \{\mathbf p^{(i)}\}_{i=1}^{N_{\text{MC}}}, \  \mathbf p^*\xrightarrow{\text{Algorithm~\ref{alg:ROB}}} \text{Output:} \ \mathbf V(\mathbf p^*)$$
\item Interpolation of PROMs:
$$\text{Input:} \ \{ \mathbf K_r (\mathbf p^{(i)})\}_{i=1}^{N_\text{MC}}, \ \{ \mathbf P^{(i,j)}\}_{i,j=1}^{N_{\text{MC}}}, \  \text{reference choice} \ i_0\xrightarrow{\text{Algorithm~\ref{alg:stepA}\&~\ref{alg:stepB}}} \text{Output:} \ \mathbf K_r(\mathbf p^*)$$
\item DMD reconstruction:
$$\text{Input:} \ \mathbf K_r(\mathbf p^*), \ \mathbf V(\mathbf p^*), \ \mathbf y(t=0; \mathbf p^*), \ \mathbf g\xrightarrow{\text{Algorithm~\ref{alg:dmd-reconst}}} \text{Output:} \ \mathbf Q_\text{DMD}(t_k;\mathbf p^*)$$
\end{itemize}
\caption{DMD-PROM framework}
\label{alg:whole}
\end{algorithm}

\begin{remark}
The sampling strategy for $\{\mathbf p^{(1)},\cdots, \mathbf p^{(N_\text{MC})}\}$ in the parameter space plays a key role in the accuracy of the subspace approximation. The so-called ``curse of dimensionality", i.e., the number of training samples $N_\text{MC}$ needed grows exponentially with the dimension of the parameter $N_\text{par}$, is a well-known challenge. In general, uniform sampling is used for $N_\text{par}\leq 5$ and moderately computationally intensive HFMs, latin hypercube sampling is used for $N_\text{par}>5$ and moderately computationally intensive HFMs, and adaptive, goal-oriented, greedy sampling is used for $N_\text{par}>5$ and highly computationally intensive HFMs. In the following numerical examples, we demonstrate our framework with $N_\text{par} =1$ for simplicity. The challenge for higher dimensions is considered as out of the scope of this work and will be studied for future work.
\end{remark}

\section{Numerical Tests}\label{sec:4}
\subsection{POD-PROM vs. DMD-PROM}
Consider a two-dimensional advection-diffusion equation
\begin{equation}\label{eq:4-1}
\frac{\partial s}{\partial t}+ \mathcal U\cdot \nabla s-\kappa \Delta s= 0,\ \mbox{for} \ \mathbf x = (x,y) \in[0,1]\times [0,1], t\in [0,1],
\end{equation}
with the boundary conditions
\begin{equation}
\begin{aligned}
&s (\mathbf x; p) =s_D(y,t), \ \mbox{for} \ \mathbf x = (x,y) \in\Gamma_D = \{0\} \times [0,1]\\
&\nabla s(\mathbf x;p)\cdot \mathbf n(\mathbf x) = 0, \ \mbox{for} \ \mathbf x\in \Gamma_N = \{\{1\}\times [0,1]\}\cup\{[0,1]\times \{1\}\}\cup \{[0,1]\times \{0\}\}
\end{aligned}
\end{equation}
where
\begin{equation}
s_D(y,t) = \left\{\begin{aligned}
&300&&\mbox{if}&y\in[0,\frac{1}{3}]\\
&300+325(\sin(3\pi |y-\bar y|)+1)&&\mbox{if}&y\in[\frac{1}{3},\frac{2}{3}]\\
&300&&\mbox{if}&y\in[\frac{2}{3},1]
\end{aligned}\right.
\end{equation}
and $\mathcal U = \begin{bmatrix}
\mathcal U_1\\
\mathcal U_2
\end{bmatrix}$, with $\mathcal U_1= p\in[0,5000]$ and $\mathcal U_2 = 0$, $\kappa=250$, $\bar y=0.4$.

The high fidelity solution is obtained via discretizing the above problem by finite differences using $75$ points in each spatial direction. Then the fully discretized variable $\mathbf S(t;p)$ is of high dimension $N_\text{sv} = 75^2$. Using upwind scheme for the advection term and implicit center difference scheme for the diffusion term with $\Delta t = 0.01$, the QoIs are chosen as $\mathbf Q(t^n;p) = \mathbf S(t^n;p)$, which follows the high fidelity linear model:
\begin{equation}\label{eq:4-1HFM}
\mathbf Q(t^{n+1};p) = \mathbf A(p)\mathbf Q(t^{n};p) +\mathbf  b(p).
\end{equation}

\begin{figure}[H]
\includegraphics[width = \textwidth]{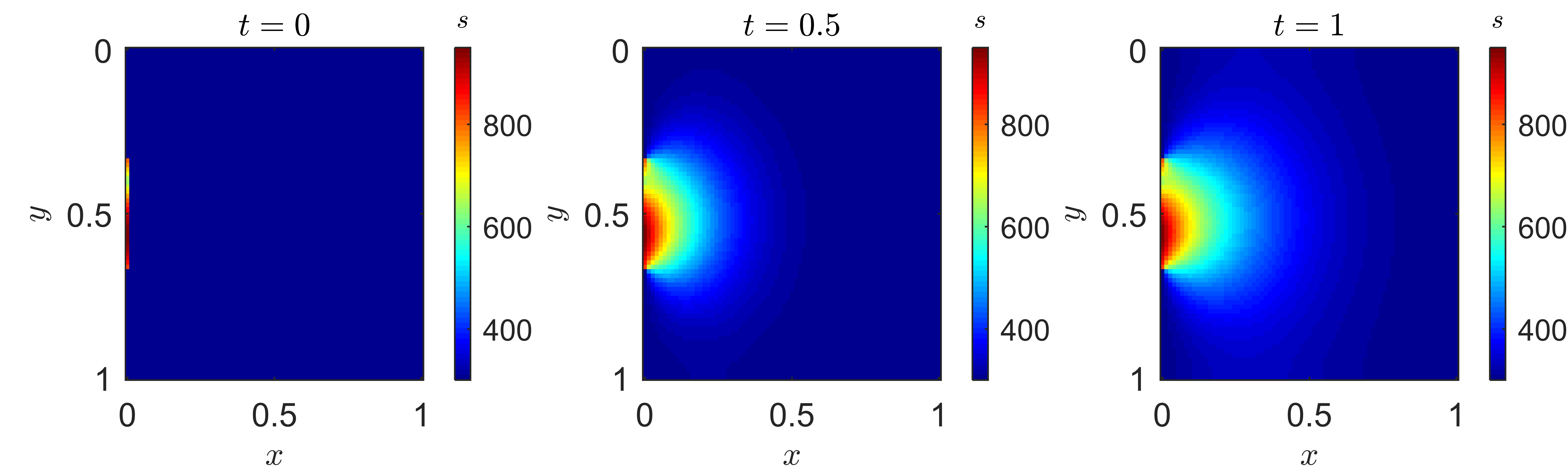}
\includegraphics[width = \textwidth]{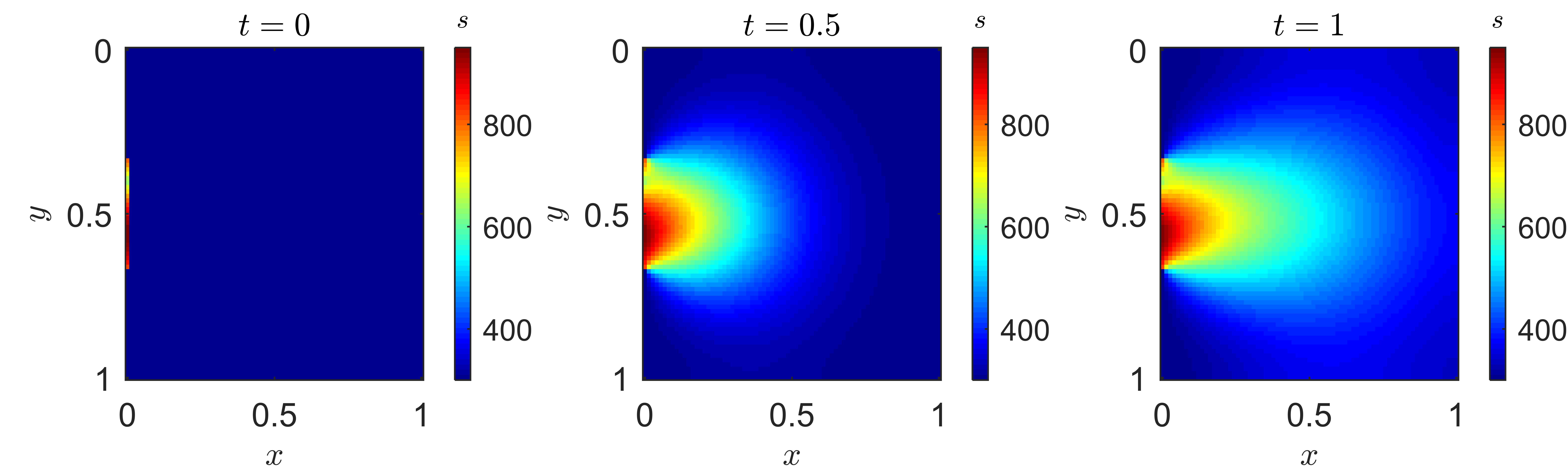}
\includegraphics[width = \textwidth]{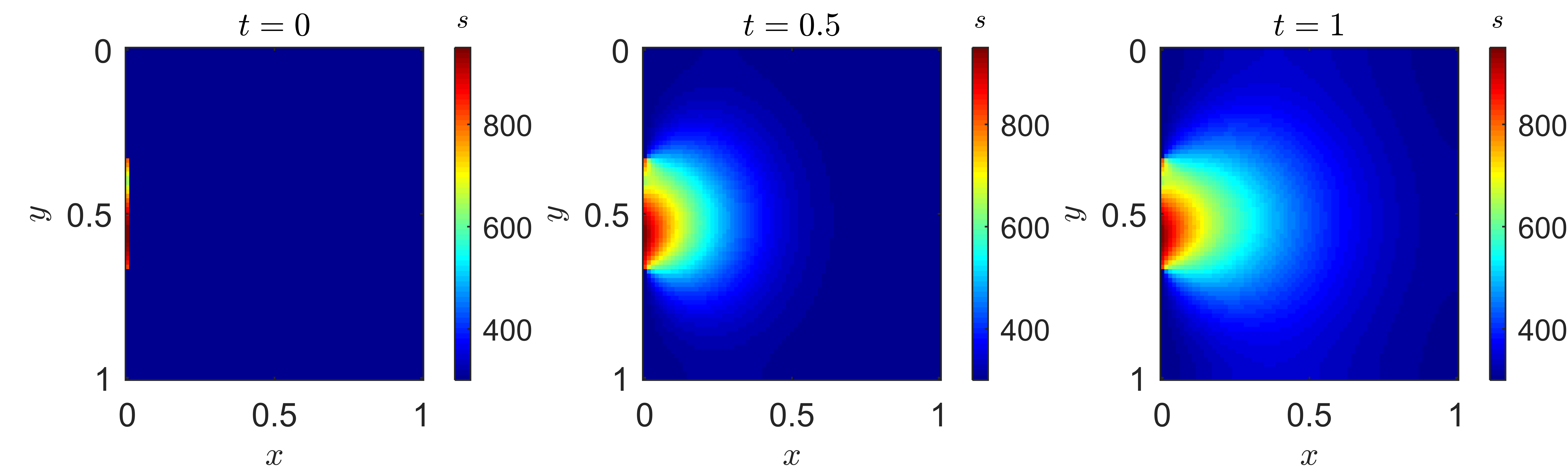}
\caption{Top row: training data from the parameter sample point $p^{(1)} = 1000$; Middle row: training data from the parameter sample point $p^{(2)} = 4000$; Bottom row: reference solution of the test data from the parameter sample point $p^{*} = 2000$. All is generated via the high fidelity simulation~\eqref{eq:4-1HFM}.}
\label{fig:4-1-1}
\end{figure}

The training data is generated from the high fidelity simulation~\eqref{eq:4-1HFM} on parameter samples $p^{(1)} = 1000$ and $p^{(2)} = 4000$. The trajectories of the solution on these two parameter samples are shown on the first two rows of Figure~\ref{fig:4-1-1}. The goal is to obtain the solution on the target parameter point $p^* = 2000$ without computing~\eqref{eq:4-1HFM}. A reference solution $\mathbf Q(t;p^*)$ on the target parameter point computed by the high fidelity simulation~\eqref{eq:4-1HFM} is shown on the bottom row of Figure~\ref{fig:4-1-1} as a reference.

We compare the performances of POD-PROM (Algorithm~\ref{alg:pod-prom}) and DMD-PROM (Algorithm~\ref{alg:whole}) in this problem. Notice that the main difference of the two methods is the way to compute the ROM $(\mathbf A_r(p^{(i)}), \mathbf b_r(p^{(i)}))$ and $\mathbf K_r(p^{(i)})$ in the offline step. In POD-PROM, $(\mathbf A_r(p^{(i)}), \mathbf b_r(p^{(i)}))$ are obtained from the Galerkin projection~\eqref{eq:s2} (or alternative projections, e.g., Petrov-Galerkin, DEIM, etc.), which requires the full knowledge of $\mathbf A(p^{(i)})$ and $\mathbf b(p^{(i)})$. However, such requirements are lifted in DMD-PROM as $\mathbf K_r(p^{(i)})$ are learned from the training dataset~\eqref{eq:data}. In this linear case, one can set the observable function $\mathbf g$ as 
\begin{equation}\label{eq:linear-obs}
\mathbf y = \mathbf g(\mathbf Q) = [1,\mathbf Q^\top]^\top,
\end{equation}
then the observable $\mathbf y$ essentially follows the following equation
\begin{equation}
\mathbf y(t^{n+1};p) = \mathbf K(p)\mathbf y(t^{n};p), \quad \mathbf K(p) = \begin{bmatrix}
1&0\\
\mathbf b(p)&\mathbf A(p)
\end{bmatrix}.
\end{equation}
Alternatively, one can use xDMD proposed in~\cite{lu2021extended} for the offline step.

Define the relative $L_2$ errors of the POD-PROM and DMD-PROM compared to the reference HFM at the test parameter point $p^*$ as 
\begin{equation}\label{eq:error1}
\begin{aligned}
&e_\text{POD}(t;p^*) = \frac{\|\mathbf Q_\text{POD}(t;p^*)-\mathbf Q(t;p^*)\|_2}{\|\mathbf Q(t;p^*)\|_2},\\
&e_\text{DMD}(t;p^*) = \frac{\|\mathbf Q_\text{DMD}(t;p^*)-\mathbf Q(t;p^*)\|_2}{\|\mathbf Q(t;p^*)\|_2}.
\end{aligned}
\end{equation}

\begin{figure}[H]
\includegraphics[width = 0.5\textwidth]{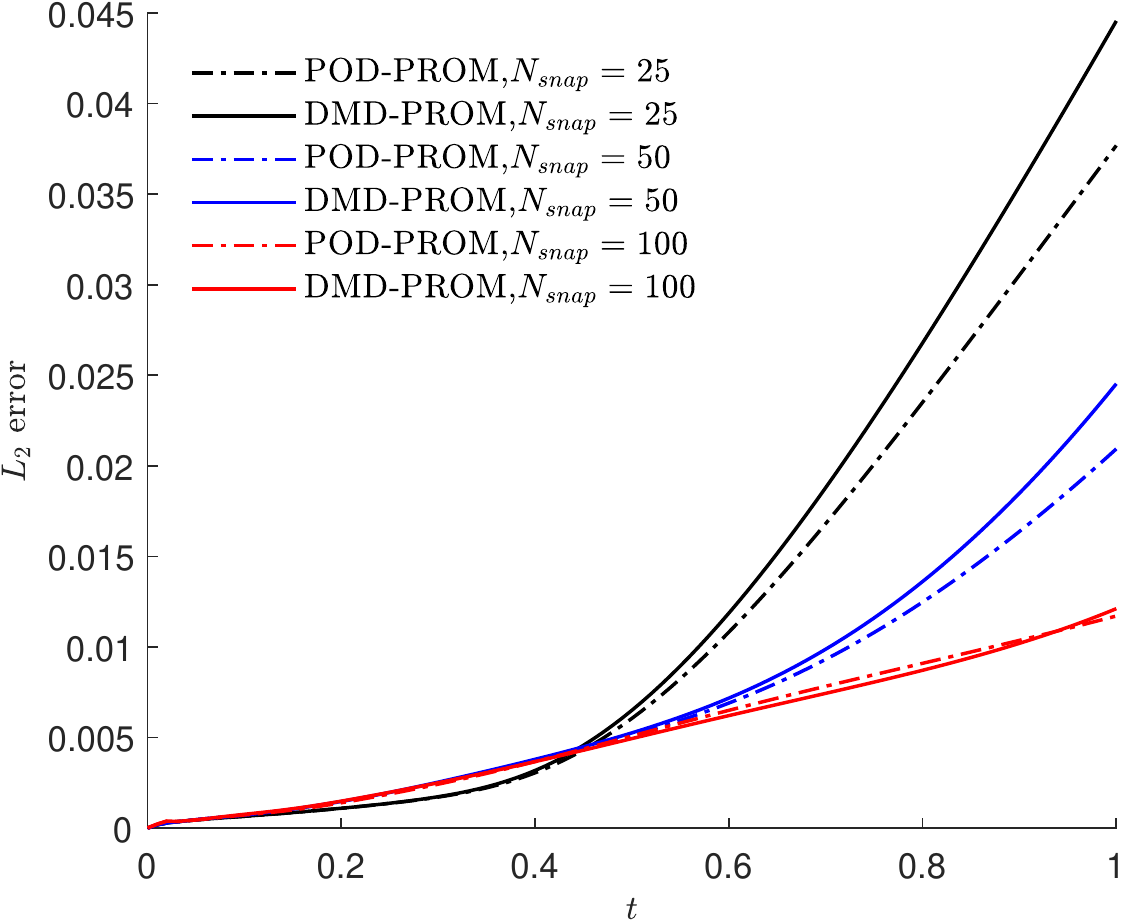}
\includegraphics[width = 0.5\textwidth]{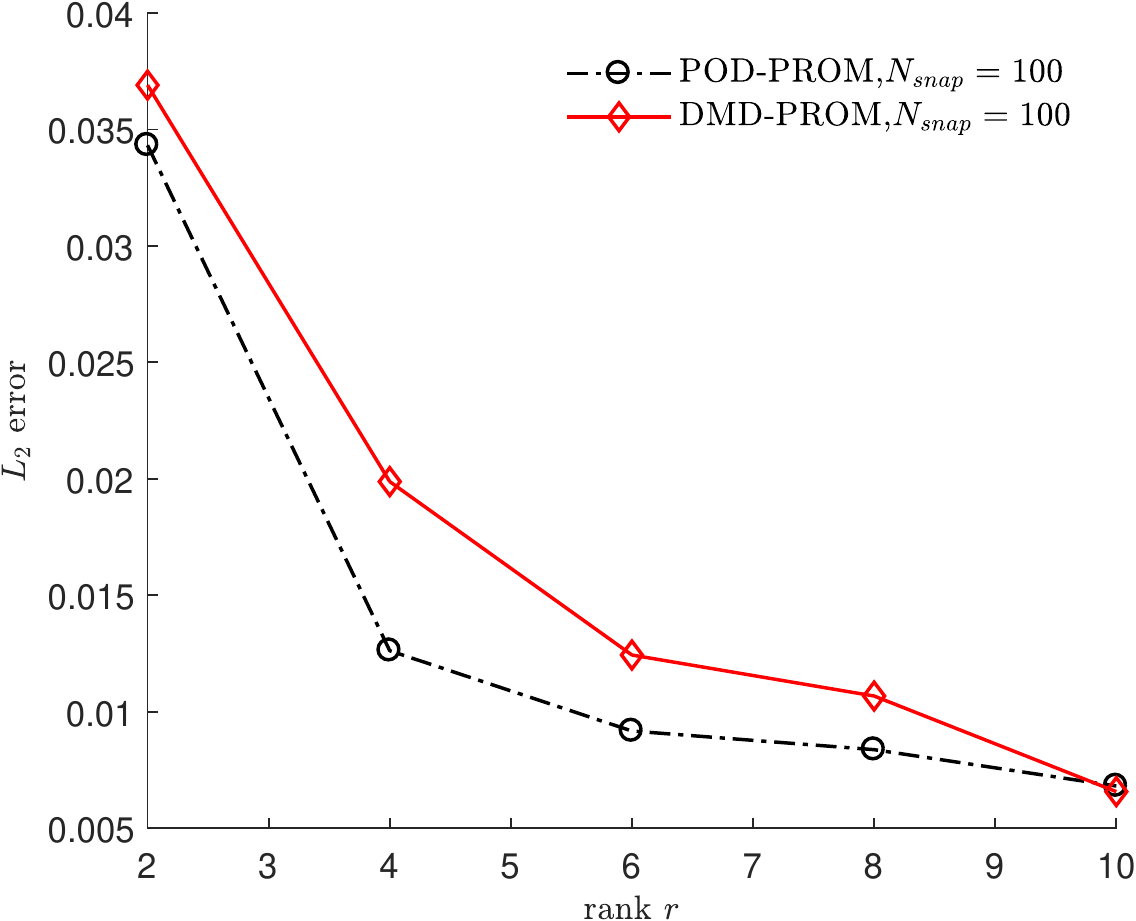}
\caption{Left: relative $L_2$ error~\eqref{eq:error1} of POD-PROMs and DMD-PROMs using different numbers of snapshots $N_\text{snap}$;  Right: total relative $L_2$ error~\eqref{eq:error2} of POD-PROMs and DMD-PROMs using different truncated rank $r$.}
\label{fig:4-1-2}
\end{figure}

We construct rank $r = 10$ POD-PROMs and DMD-PROMs using different numbers of snapshots $N_\text{snap}$. The relative $L_2$ error of each are displayed in the left plots of Figure~\ref{fig:4-1-2}. For both methods, the magnitudes of the errors decrease with larger number of $N_\text{snap}$. In POD, this can be explained by the fact that the low dimensional subspace (represented by the ROBs determined by the training data) becomes more optimal with larger number of informative data. The same reason holds for DMD. Furthermore, one can observe that DMD has larger error than POD in the case $N_\text{snap} = 25$ and $N_\text{snap} =50$. The difference is caused by the loss of accuracy in learning $\mathbf K_r(p^{(i)})$ from data in the DMD framework. Given enough data (in the case $N_\text{snap} = 100$), the difference between POD and DMD saturates as the learned $\mathbf K_r(p^{(i)})$ becomes accurate enough.

Define the total relative $L_2$ errors of the POD-PROM and DMD-PROM as
\begin{equation}\label{eq:error2}
\begin{aligned}
&E_\text{POD}(p^*) = \frac{\sqrt{\sum_{k=1}^{N_T}\|\mathbf Q_\text{POD}(t_k;p^*)-\mathbf Q(t_k;p^*)\|_2^2}}{\sqrt{\sum_{k=1}^{N_T}\|\mathbf Q(t_k;p^*)\|_2^2}},\\
&E_\text{DMD}(p^*) = \frac{\sqrt{\sum_{k=1}^{N_T}\|\mathbf Q_\text{DMD}(t_k;p^*)-\mathbf Q(t_k;p^*)\|_2^2}}{\sqrt{\sum_{k=1}^{N_T}\|\mathbf Q(t_k;p^*)\|_2^2}}.
\end{aligned}
\end{equation}

We construct POD-PROMs and DMD-PROMs of different ranks $r$ using the same number of snapshots $N_\text{snap}=N_T = 100$. The total relative $L_2$ error of each are displayed in the right plots of Figure~\ref{fig:4-1-2}. Both methods have decreased errors with larger rank $r$ and the errors saturate after $r = 10$, which verifies the low rank nature of this problem. DMD has slightly larger error than POD for small ranks because of the loss in accuracy of learning $\mathbf K_r(p^{(i)})$. The accuracy of the DMD-based ROM is affected by the number of snapshots and the rank truncation, as studied by many previous papers, e.g., ~\cite{lu2020prediction}.

The study of this example illustrates that the knowledge of the physics (i.e., the specific value of $\mathbf A(p^{(i)})$ and $\mathbf b(p^{(i)})$) is a privilege in constructing more accurate ROMs (as conventionally done in POD). However, one can still construct satisfactory ROMs without full knowledge of the physics as long as enough data is given. Our proposed DMD-PROM framework compensates the short of knowledge of the physics by a data-driven learning. The data-driven learning part can be improved by partial knowledge of the physics via the observable function $\mathbf g$, which has been demonstrated in many previous studies~\cite{lu2020prediction,lu2020lagrangian,lu2021dynamic,lu2021extended,nathan2018applied,williams2015data}. The following numerical examples will show how to incorporate the physics ingredients into the framework.

\subsection{Linear Diffusion Equation} 
\begin{figure}[H]
\begin{center}
\includegraphics[trim={0 0.1cm 6cm 0},clip]{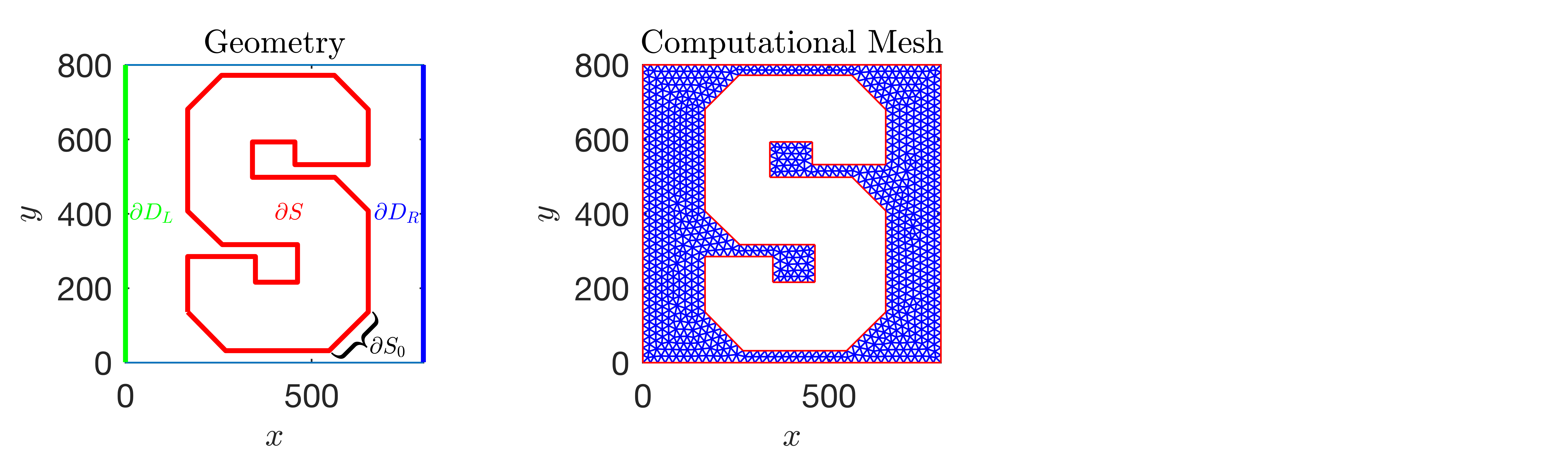}
\end{center}
\caption{Multi-connected simulation domain $\mathcal D$ (left) and the mesh used in the finite-element solution of~\eqref{eq:test2}.}
\label{fig4}
\end{figure}

Consider a two-dimensional diffusion equation in a multi-connected domain $\mathcal D$ with inhomogeneous boundary conditions,
\begin{equation}\label{eq:test2}
\left\{
\begin{aligned}
&\rho c\frac{\partial s}{\partial t} -\nabla\cdot(k\nabla s)=0,  \qquad (x,y) \in \mathcal D, \quad t\in (0,5000];\quad \rho = 1, k = 1, c = p\in [1,2];\\
 &s(x,y,0) = 0;\\
 & s(0,y,t) =  2, \quad s(800,y,t) = 1, \quad \frac{\partial s}{\partial y}(x,0,t) = \frac{\partial s}{\partial y}(x,800,t) = 0, \quad 
  s(x,y,t) =  3 \;\;\;  \mbox{on $\partial \mathcal S$ (red)}.\\
 \end{aligned}\right.
 \end{equation}
The domain $\mathcal D$ is the $800 \times 800$ square with an S-shaped cavity (Figure~\ref{fig4}). The Dirichlet boundary conditions are imposed on the left and right sides of the square and the cavity surface. The top and bottom of the square are impermeable.  The training data and reference solutions are obtained via  Matlab PDE toolbox on the finite-element mesh with $4287$ elements shown in Figure~\ref{fig4}. 

\begin{figure}[H]
\includegraphics[width = \textwidth]{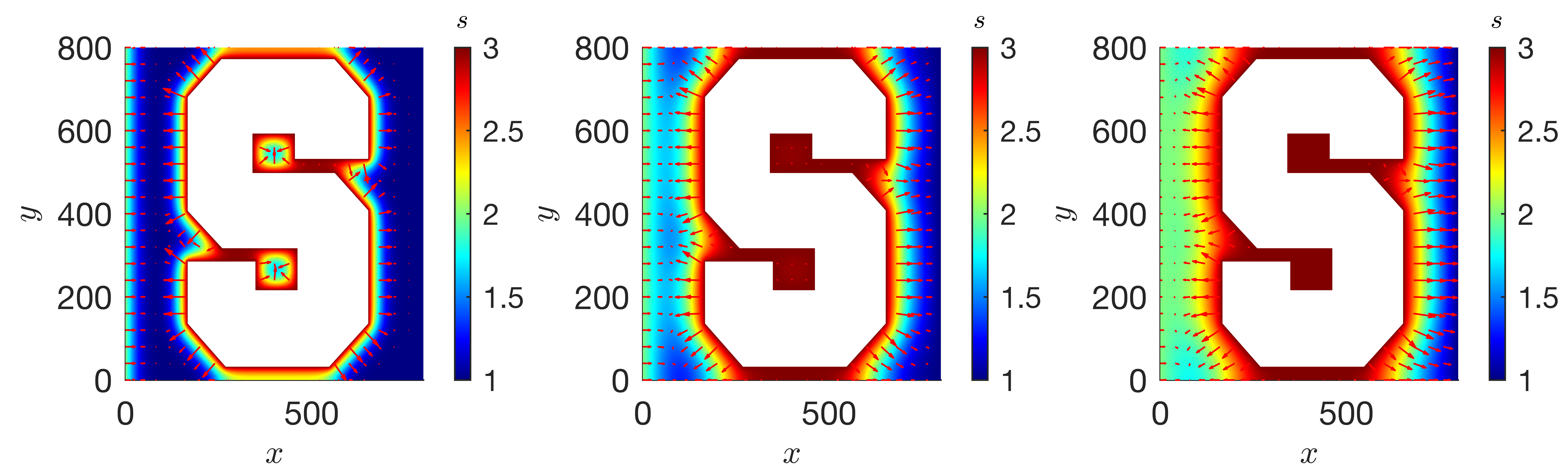}
\includegraphics[width = \textwidth]{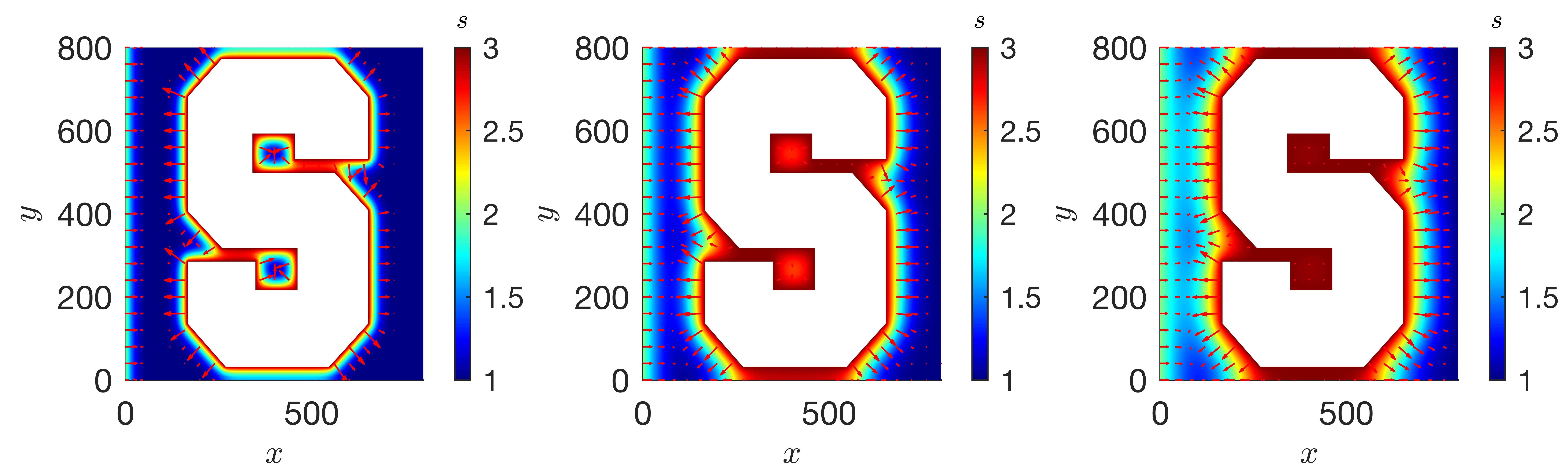}
\includegraphics[width = \textwidth]{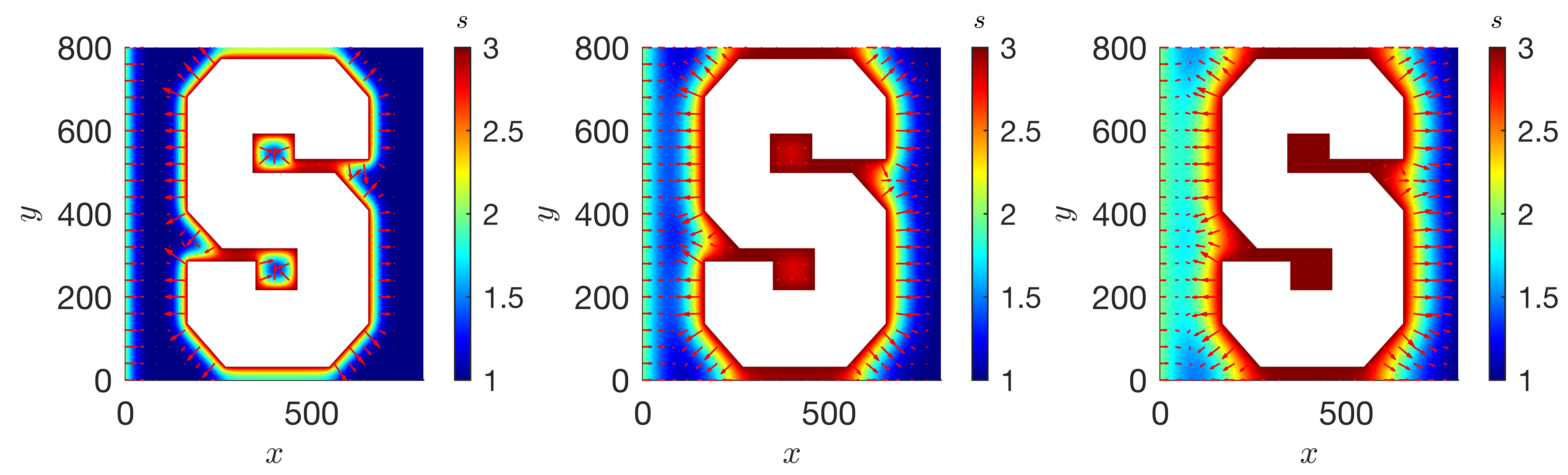}
\includegraphics[width = \textwidth]{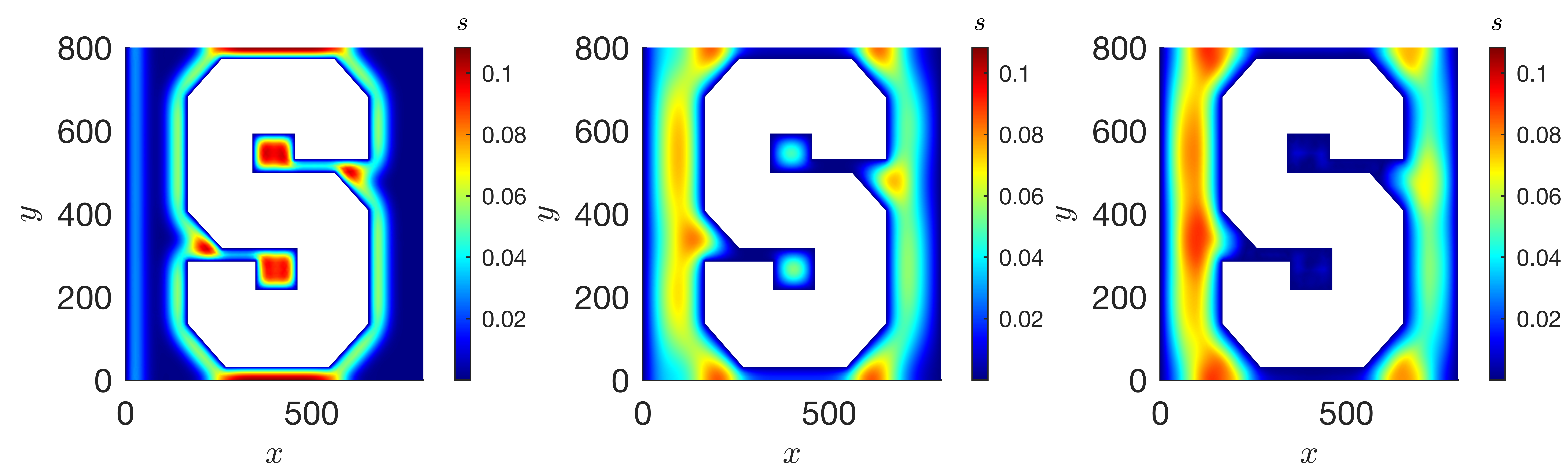}
\caption{The training data at $p^{(1)} = 1$ (first row), $p^{(2)} = 2$ (second row), the reference solution at target parameter point $p^* = 1.5$ (third row) and its corresponding DMD-PROM absolute error map (last row) are shown at times $t = 500$ (left), $t = 2500$ (middle) and $t=5000$(right) for problem~\eqref{eq:test2}.}
\label{fig5}
\end{figure}

\begin{figure}[H]
\includegraphics[width = 0.33\textwidth]{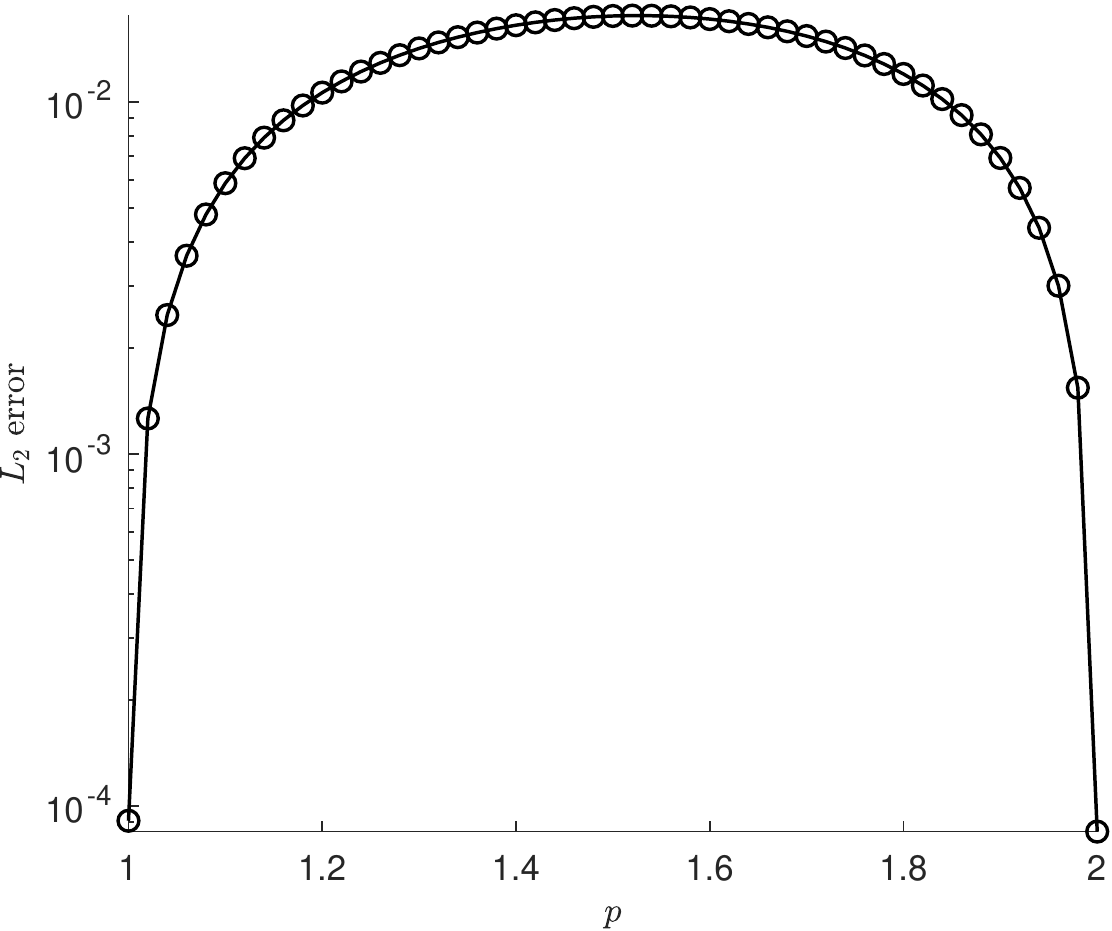}
\includegraphics[width = 0.33\textwidth]{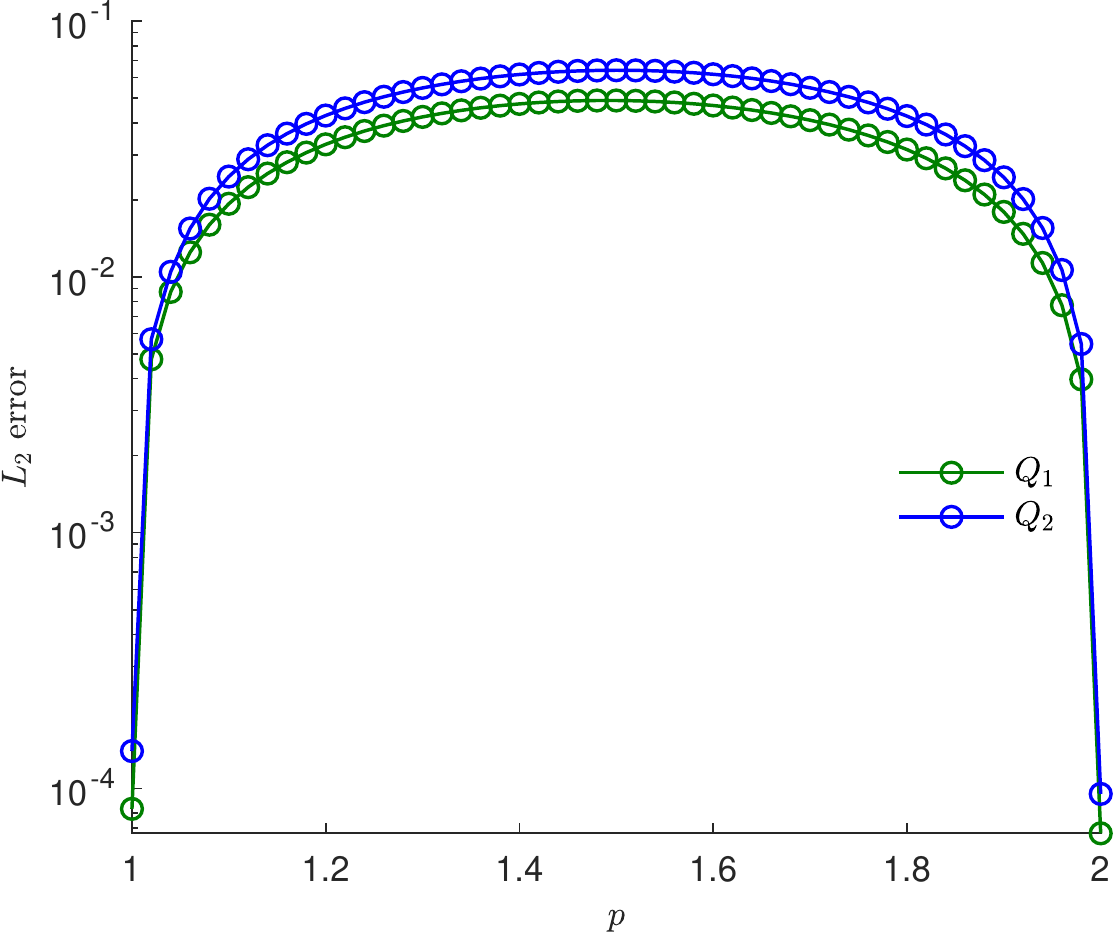}
\includegraphics[width = 0.33\textwidth]{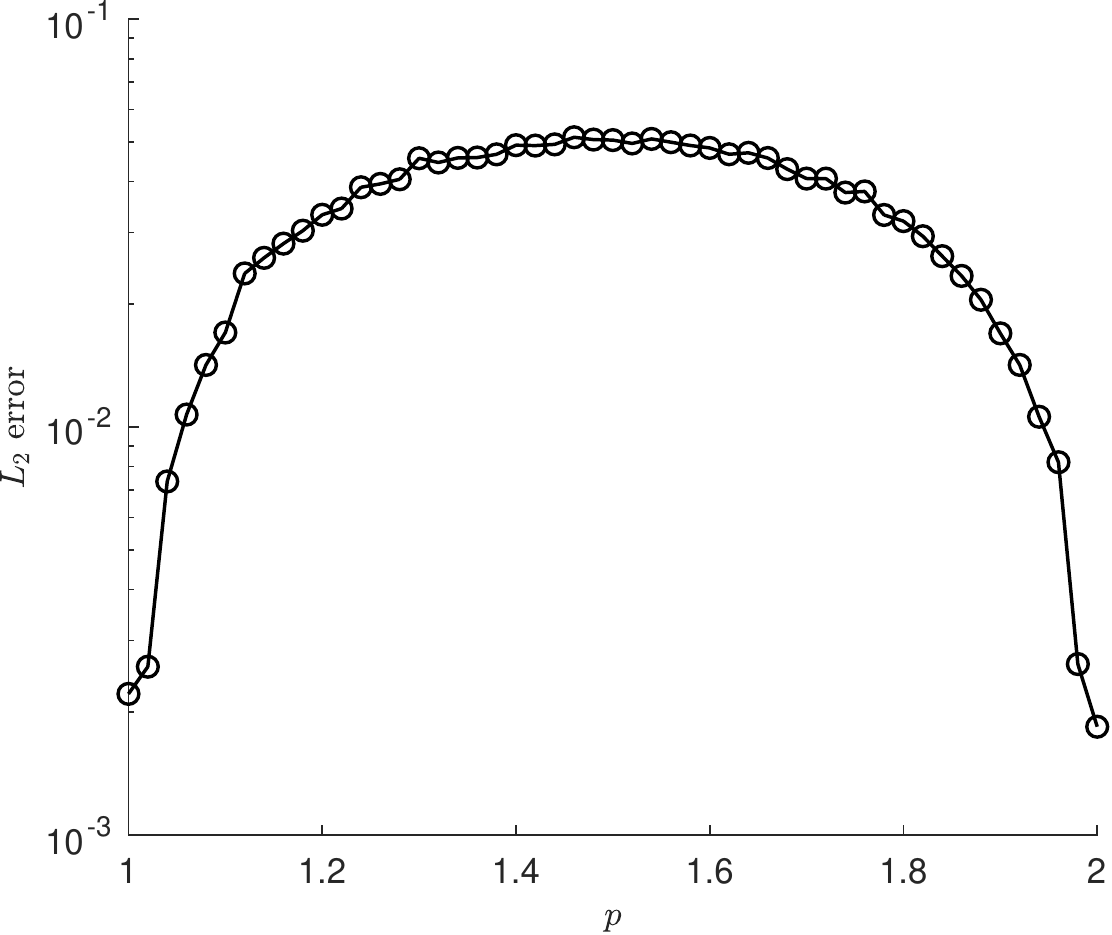}
\caption{Left: Total $L_2$ error~\eqref{eq:error2} of $r=10$ rank DMD-PROM for QoIs $\mathbf Q = \mathbf S$ with respect to different values of $p^*\in [1,2]$;  Middle: Total $L_2$ error~\eqref{eq:error2} of $r=20$ rank DMD-PROM for QoIs $\mathbf Q_1 = -k\partial_x\mathbf S$ and $\mathbf Q_2 = -k\partial_y\mathbf S$ with respect to different values of $p^*\in [1,2]$; Right: Total $L_2$ error~\eqref{eq:error2} of $r=8$ rank DMD-PROM for the heat rate across the edge of interest $\partial S_0$ with respect to different values of $p^*\in [1,2]$.}
\label{fig7}
\end{figure}

First, we consider the case where the QoIs are chosen as $\mathbf Q(t^n;p) = \mathbf S(t^n;p)$. The observable function $\mathbf g$ is chosen the same as~\eqref{eq:linear-obs}. The training data is generated on the parameter samples $p^{(1)} = 1$ and $p^{(2)} = 2$ and the reference solution of a particular test parameter point $p^* = 1.5$, together with the absolute error map of the DMD-PROM, are presented in Figure~\ref{fig5}. Without the explicit formula of the HFM, which is designed in the Matlab PDE toolbox, our framework can provide an accurate ROM with small rank $r = 10$. The Matlab PDE toolbox essentially generates the numerical solutions by solving a HFM iteration in a form of~\eqref{eq:4-1HFM}. If POD framework is used instead, one has to know the exact values of $\mathbf A(p)$ and $\mathbf b(p)$, which require large numerical differential matrices and embeddings of the complex boundary conditions. The left plot of Figure~\ref{fig7} shows the accuracy of the $10$ rank DMD-PROM for different values of $p^*\in [1,2]$. The error plot has a peak in the middle of the parameter space and decays to its minimum at the two ends, where the errors only come from the differences between the HFMs and the ROMs constructed by local ROBs (at $p^{(1)} = 1$ and $p^{(2)} = 2$ respectively).

Then we consider the case where the QoIs are chosen as the heat flux, i.e.,  $\mathbf Q_1(t^n;p) = -k\partial_x\mathbf S(t^n;p)$ and $\mathbf Q_2(t^n;p) = -k\partial_y\mathbf S(t^n;p)$. The heat flux vector fields of the training data and test reference solutions are plotted using arrows in Figure~\ref{fig5}. The observable function $\mathbf g$ is chosen as 
\begin{equation}
\mathbf y = \mathbf g(\mathbf Q_1,\mathbf Q_2) = [1,\mathbf Q_1^\top,\mathbf Q_2^\top]^\top.
\end{equation}
The middle plot of Figure~\ref{fig7} shows the accuracy of the $20$ rank DMD-PROM for different values of $p^*\in [1,2]$.  Similar error distributions are observed as before.

Finally, we consider the case where the QoIs are chosen as the heat rate across the edge of interest (shown in Figure~\ref{fig4}), i.e., $Q(t^n;p) = \int_{\partial S_0}-k\nabla s(t^n,\mathbf x;p)\cdot \vec{\mathbf n}_{\partial S_0}dA$. The observable function $\mathbf g$ is chosen as Hermite polynomials $\mathcal H$ of $Q(t^n;p)$ up to order $m = 8$,
\begin{equation}
\mathbf y = \mathbf g(Q) = [\mathcal H_0(Q),\cdots, \mathcal H_{m-1}(Q)]^\top.
\end{equation}
The right plot of Figure~\ref{fig7} shows the accuracy of the $8$ rank DMD-PROM for different values of $p^*\in [1,2]$.  Similar error distributions are observed as before. The data-driven learning ability of the DMD allows one to construct the ROM for the QoIs directly without accessing the expensive HFM and post-processing the high dimensional state variables. 

\subsection{Nonlinear Navier-Stokes Equation}\label{sec:NS}
We consider two-dimensional flow of an incompressible fluid with density $\rho = 1$ and dynamic viscosity $\nu =  p \in [1/400,1/800]$ (these and other quantities are reported in consistent units) around an impermeable circle of diameter $D = 0.1$. The flow, which takes place inside a rectangular domain $\mathcal D = \{ \mathbf x = (x,y)^\top : (x,y) \in [0,2] \times [0,1] \}$, is driven by an externally imposed pressure gradient; the center of the circular inclusion is $\mathbf x_\text{circ} = (0.3,0.5)^\top$. Dynamics of the three state variables, flow velocity $\mathbf s(\mathbf x,t) = (s_1,s_2)^\top$ and fluid pressure $P(\mathbf x,t)$, is described by the two-dimensional Navier-Stokes equations,
\begin{equation}\label{eq:NS}
\left\{
\begin{aligned}
&\frac{\partial s_1}{\partial x}+\frac{\partial s_2}{\partial y} = 0; \\
&\frac{\partial s_1}{\partial t}+s_1\frac{\partial s_1}{\partial x}+s_2\frac{\partial s_1}{\partial y} = -\frac{1}{\rho}\frac{ \partial P}{ \partial x}+\nu\left(\frac{\partial^2 s_1}{\partial x^2}+\frac{\partial^2 s_1}{\partial y^2}\right), \qquad \mathbf x \in \mathcal D, \quad t > 0; \\
&\frac{\partial s_2}{\partial t}+u\frac{\partial s_2}{\partial x}+s_2\frac{\partial s_2}{\partial y} = -\frac{1}{\rho}\frac{\partial P}{\partial y}+\nu\left(\frac{\partial^2 s_2}{\partial x^2}+\frac{\partial^2 s_2}{\partial y^2}\right);
\end{aligned}
\right.
\end{equation}
subject to initial conditions $\mathbf s(x,y,0) = (0,0)^\top$ and $P(x,y,0) = 0$; and boundary conditions 
\[
P(2,y,t) = 0, \quad \frac{\partial P}{\partial \mathbf n}|_{\partial \mathcal D\setminus \{x= 2\}} = 0, \quad 
\mathbf s(0,y,t) = (1,0)^\top, \quad \frac{\partial \mathbf s(2,y,t)}{\partial \mathbf n} = 0, \quad \mathbf s(x,0,t) = \mathbf s(x,1,t) = \mathbf 0.
\] 
Here $\mathbf n$ denotes the unit normal vector. 

The training data and reference solution are obtained with the Matlab code~\cite{NScode}, which implements a finite-difference scheme on the staggered grid with $\Delta x = \Delta y = 0.02$ and $\Delta t = 0.0015$. Our QoIs are the magnitude of the flow velocity, $\mathbf Q(t^n;p) = \sqrt{\mathbf S_1(t^n;p)^2+\mathbf S_2(t^n;p)^2}$, which follows a nonlinear equation in an implicit form. We collect $N_\text{snap} = 250$ snapshots of $\mathbf Q$ from $t = 4.125$ to $t = 4.5$ operated at $p^{(1)} = 1/400$ and $p^{(2)} = 1/800$ into the training dataset (shown in Figure), from which DMD learns the nonlinear dynamics from a $10$-rank surrogate using the same observable as~\eqref{eq:linear-obs}. This construction is equivalent to the xDMD and the same numerical example was shown in section 4.3 of~\cite{lu2021extended}. The reference solution of a particular test parameter point $p^* =1.5$, together with the absolute error map of the DMD-PROM, are presented in Figure~\ref{fig9}.

Figure~\ref{fig10} shows the accuracy of the $10$-rank DMD-PROM for different values of $p^*\in[1/400,1/800]$. Another alternative approach to tackle this problem is Kriging. Kriging interpolates the velocity magnitude map pixel by pixel based on prior covariances, which depend on the distance between the target point and the sampled points in the parameter space. It is a pure data-driven approach in the sense that it does not take the underlying governing equations into account. Although Kriging is widely used in many applications (most linear problems) successfully, previous study~\cite{amsallem2011online} demonstrates that such pure data-driven approach may fail to detect the bifurcations in a complex dynamic system. Figure~\ref{fig9} and~\ref{fig10} also verify the disadvantage of Kriging in accurately capturing the complex nonlinear dynamics compared to the physics-aware DMD-PROM.

\begin{figure}[H]
\includegraphics[width = \textwidth]{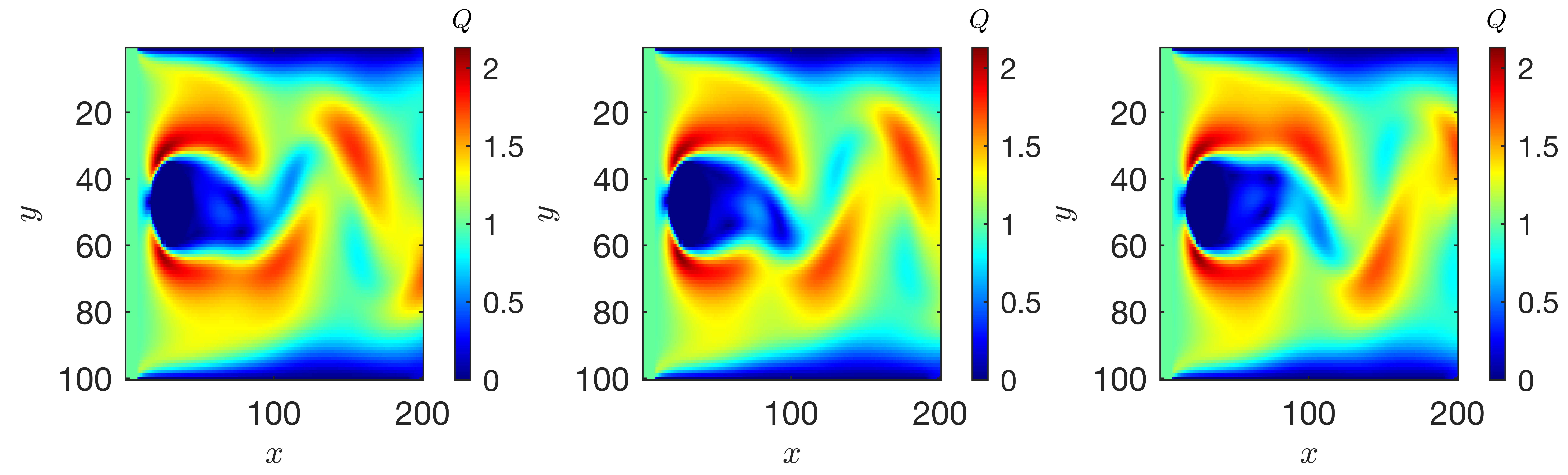}
\includegraphics[width = \textwidth]{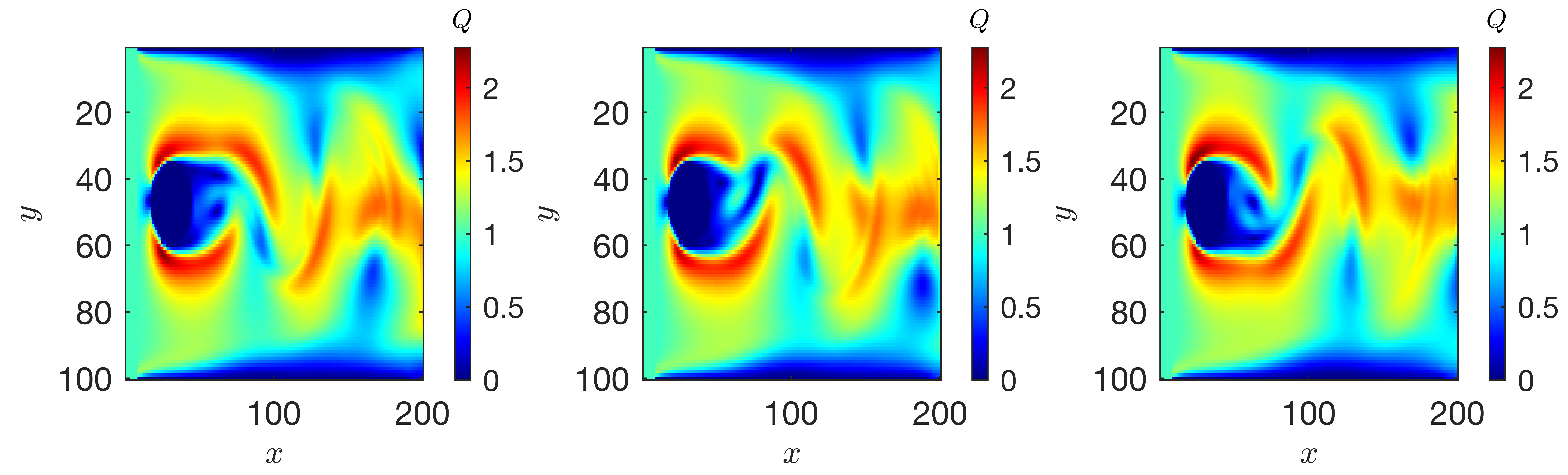}
\caption{The training data of~\eqref{eq:NS} at $p^{(1)} = 1/400$ (top row) and $p^{(2)} = 1/800$ (bottom row), at times $t = 4.125$ (left), $t =4.3125$ (middle) and $t=4.5$(right).}
\label{fig8}
\end{figure}

\begin{figure}[H]
\includegraphics[width = \textwidth]{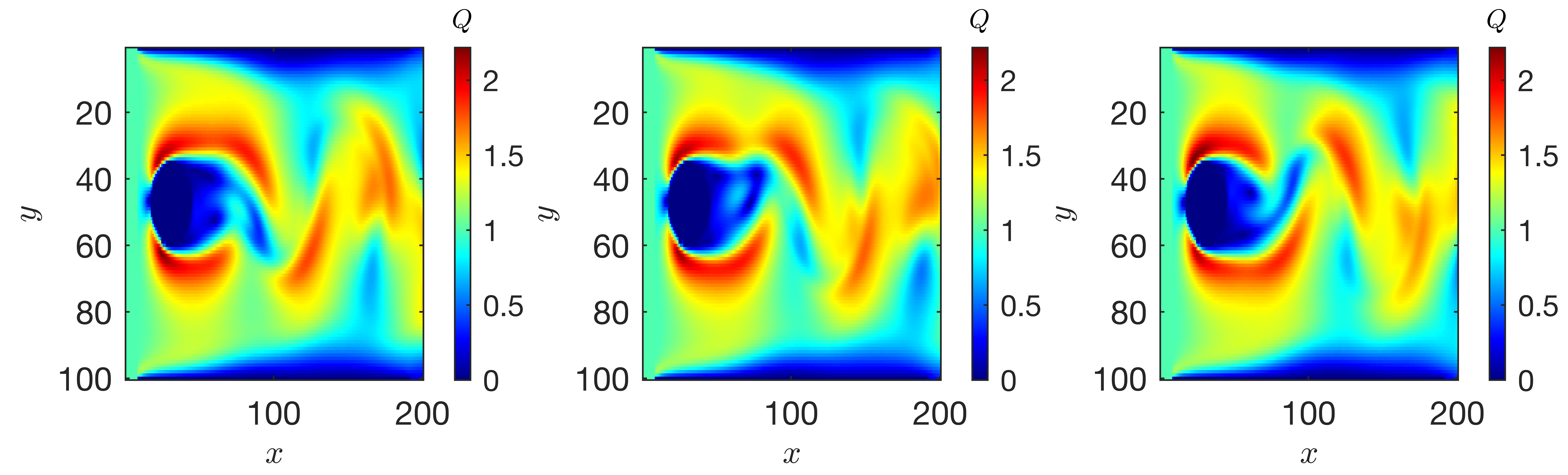}
\includegraphics[width = \textwidth]{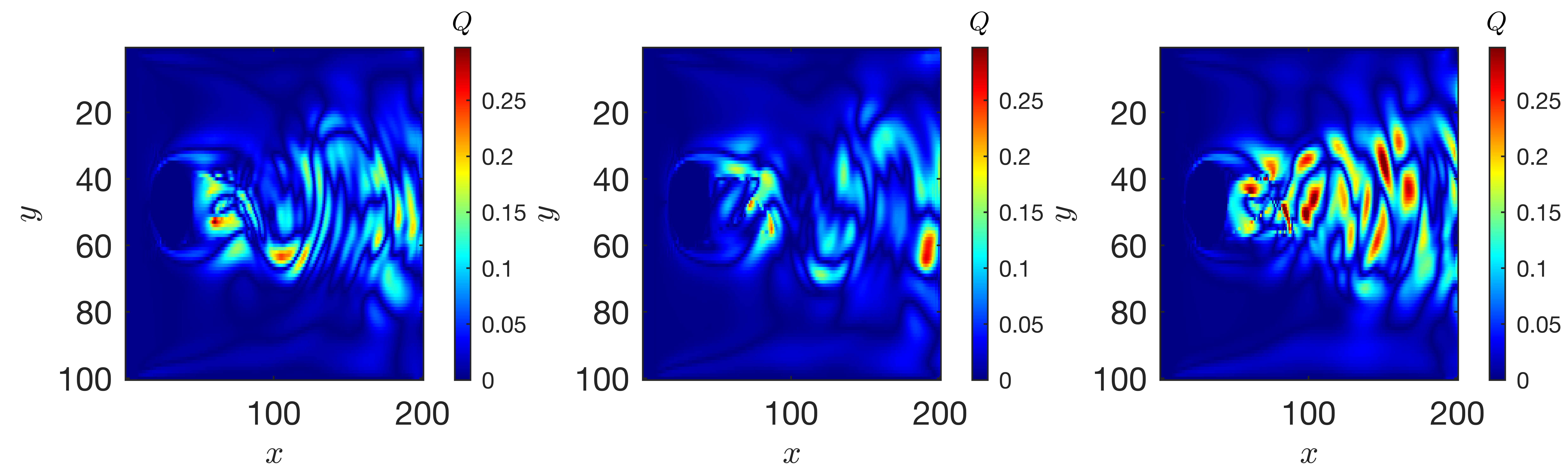}
\includegraphics[width = \textwidth]{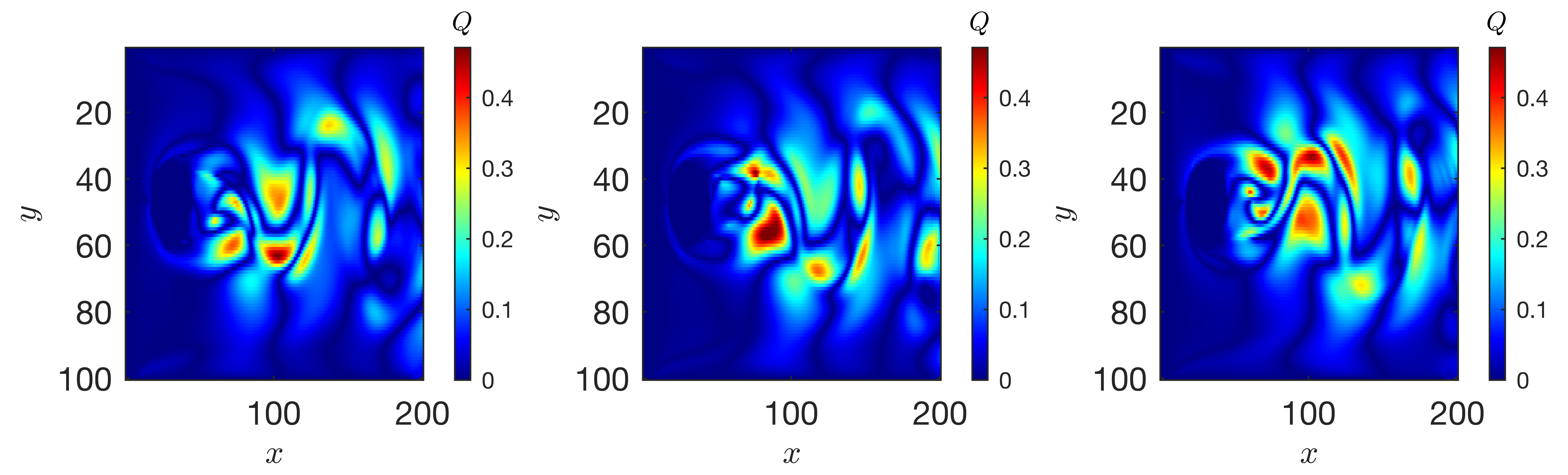}
\caption{The reference solution (top row) of~\eqref{eq:NS}, corresponding DMD-PROM absolute error map (middle row)  and Kriging absolute error map (bottom row) at target parameter point $p^* = 1/600$, at times $t = 4.125$ (left), $t = 4.3125$ (middle) and $t=4.5$ (right).}
\label{fig9}
\end{figure}

Table~\ref{table1} shows the comparison of computational costs in exploring the dynamics with different parameter values $p^*\in [1/800,1/400]$ using the HFM, Kriging and DMD-PROM respectively. The parameter space of interests is discretized into $51$ uniformly distributed points, including the $2$ endpoints (used as training parameter samples) and the rest $49$ inner points (used as test parameter points). Although Kriging has a larger online speedup, the proposed DMD-PROM has the same level total speedup as Kriging and much better accuracy in all test parameter points as shown in Figure~\ref{fig10}. 

\begin{minipage}[t]{.4\textwidth}
\begin{figure}[H]
\includegraphics[width = \textwidth]{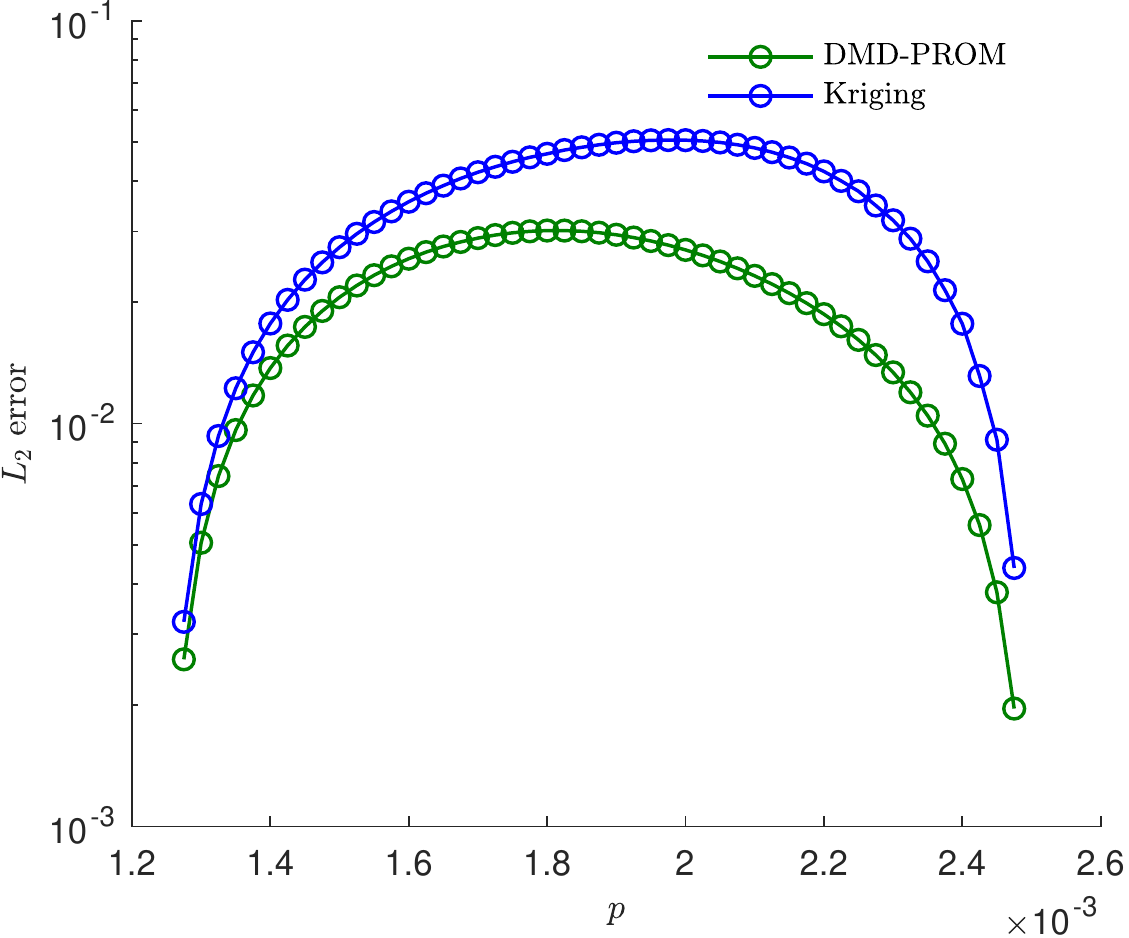}
\caption{Total $L_2$ error~\eqref{eq:error2} of $r= 10$ rank DMD-PROM (blue) and Kriging (green) for QoIs $\mathbf Q = \sqrt{\mathbf S_1(t^n;p)^2+\mathbf S_2(t^n;p)^2}$ with respect to different values of $p^*\in [1/800,1/400]$.}
\label{fig10}
\end{figure}
\end{minipage}
\hspace{0.2cm}
\begin{minipage}[t]{.55\textwidth}
\begin{center}
\vspace{2cm}
\captionof{table}{Comparison of computational performances}
\label{table1}
\begin{tabular}{ |c|c|c|c|}
 \hline
 Approach& HFM& Kriging& DMD-PROM \\ 
 \hline
 Offline(seconds)&  --&$15.2195$&$14.0906$\\ 
 \hline
Online (seconds) &$255.1326$ &$1.5288$ &$2.8567$\\ 
 \hline
Online speedup&--&$167$ &$89$\\
 \hline
 Total speedup&--& $15$&$15$\\
 \hline
\end{tabular}
\end{center}
\end{minipage}

\section{Conclusion}\label{sec:5}
In this work, we propose a new physics-aware data-driven framework of model reduction for parametric complex systems. It combines the advantages of the popular data-driven modeling tool DMD and previous studies on manifold interpolations over the parameter space of interests. The novelty of this work includes: 1. alleviating the strict constraints on the access to the HFMs of the high dimensional state variables in projection-based PROM methods by constructing low-dimensional surrogate models for the QoIs directly in a physics-aware data-driven way; 2. improving the PROM accuracy from conventional pure-data driven tools like Kriging and requiring much less training data from modern nonlinear ML tools like deep learning. The proposed framework shows robustness and flexibility by providing a bridge between the understanding of data and physics. The significant speedup also verifies that DMD-PROM is suitable for real-time processing. 

In the follow-up studies we will develop the current framework in the following directions: 1). To further improve the online speedup, more efficient manifold interpolation methods (e.g., ~\cite{zhang2021gaussian}) need to be explored; 2). For problems in high-dimensional parameter spaces, adaptive sampling strategies need to be carefully designed and the current framework needs to be applied together with practical considerations of parameter space reduction~\cite{lieberman2010parameter}; 3). Further studies on error estimation of the whole framework are needed in order to construct reliable surrogate models, which may further enable outer-loop applications such as design, inverse problems, optimization and uncertainty quantification. 

\newpage

\begin{center}
\textbf{SUPPLEMENTAL MATERIAL}
\end{center}
\setcounter{equation}{0}
\setcounter{figure}{0}
\setcounter{table}{0}
\setcounter{algocf}{0}
\makeatletter
\renewcommand{\theequation}{S\arabic{equation}}
\renewcommand{\thefigure}{S\arabic{figure}}
\renewcommand{\thealgocf}{S\arabic{algocf}}

\begin{algorithm}[H]
\textit{Offline Step:} 

\textbf{For} $i = 1,\cdots,N_\text{MC}$,

\ \ \ \ \ \ Compute the high fidelity training data~\eqref{eq:data},

\ \ \ \ \ \ Input: $\{\mathbf Q(t_0;\mathbf p^{(i)}),\cdots,\mathbf Q(t_{N_\text{snap}};\mathbf p^{(i)})\}$, $\mathbf A(p^{(i)})$ and $\mathbf b(p^{(i)})$, 
\begin{enumerate}
\item Apply SVD $[\mathbf Q(t_0;\mathbf p^{(i)})^\top,\cdots,\mathbf Q(t_{N_\text{snap}};\mathbf p^{(i)})^\top]^\top\approx \mathbf V(\mathbf p^{(i)}) \boldsymbol \Sigma(\mathbf p^{(i)}) \mathbf Z (\mathbf p^{(i)}) ^*$ with $\mathbf V(\mathbf p^{(i)}) \in \mathbb C^{N\times r}$, $\boldsymbol \Sigma(\mathbf p^{(i)})  \in \mathbb C^{r\times r}$, $\mathbf Z(\mathbf p^{(i)}) \in \mathbb C^{r\times N_\text{snap}}$, where $r$ is the truncated rank chosen by certain criteria and should be the same for all $i = 1,\cdots, N_\text{MC}$.
\item Use Galerkin-projection to compute local ROMs:
\begin{equation}\label{eq:s2}
\mathbf A_r(p^{(i)}) = \mathbf V(p^{(i)})^\top\mathbf A(p^{(i)}) \mathbf V(p^{(i)}), \mathbf b_r(p^{(i)}) = \mathbf V(p^{(i)})^\top\mathbf b(p^{(i)}),
\end{equation}
\item Compute $\mathbf P^{(i,j)} = \mathbf V(\mathbf p^{(i)})^\top\mathbf V(\mathbf p^{(j)})$ for $i,j = 1,\cdots, N_\text{MC}$.
\end{enumerate}
\ \ \ \ \ \ Output: $\mathbf V (\mathbf p^{(i)})$, $\mathbf P^{(i,j)}$, $\mathbf A_r(\mathbf p^{(i)})$ and $\mathbf b_r(\mathbf p^{(i)})$.

\textbf{End}

\textit{Online Step:}
\begin{itemize}
\item Interpolation of ROBs:
$$\text{Input:} \ \{ \mathbf V (\mathbf p^{(i)})\}_{i=1}^{N_{\text{MC}}}, \ \{ \mathbf P^{(i,j)}\}_{i,j=1}^{N_{\text{MC}}}, \ \{\mathbf p^{(i)}\}_{i=1}^{N_{\text{MC}}}, \  \mathbf p^*\xrightarrow{\text{Algorithm~\ref{alg:ROB}}} \text{Output:} \ \mathbf V(\mathbf p^*)$$
\item Interpolation of PROMs:
$$\text{Input:} \ \{ \mathbf A_r (\mathbf p^{(i)}),\mathbf P^{(i,j)}\}_{i,j=1}^{N_{\text{MC}}}, \  \text{reference choice} \ i_0\xrightarrow{\text{Algorithm~\ref{alg:stepA}\&~\ref{alg:stepB}}} \text{Output:} \ \mathbf A_r(\mathbf p^*)$$
$$\text{Input:} \ \{ \mathbf b_r (\mathbf p^{(i)}), \mathbf P^{(i,j)}\}_{i,j=1}^{N_{\text{MC}}}, \  \text{reference choice} \ i_0\xrightarrow{\text{Algorithm~\ref{alg:stepA}\&~\ref{alg:stepB}}} \text{Output:} \ \mathbf b_r(\mathbf p^*)$$
\item POD reconstruction:

Input: $\mathbf A_r(\mathbf p^{*})$, $\mathbf b_r(\mathbf p^{*})$, $\mathbf V(\mathbf p^{*})$ and $\mathbf Q(t_0;\mathbf p^{*})$,

\ \ \ \ \ \ \ \ \ \ $\mathbf q(t_0;\mathbf p^*) = \mathbf V(\mathbf p^*)^\top\mathbf Q(t_0;\mathbf p^{*})$

\ \ \ \ \ \ \ \ \ \ \textbf{For} $k = 1,\cdots,N_T$,

\ \ \ \ \ \ \ \ \ \ \ \ \ \ \ \ $\mathbf q(t_{k};\mathbf p^*) = \mathbf A_r(p^{*})\mathbf q(t_{k-1};\mathbf p^*) + \mathbf b_r(p^{*}),\ \mathbf Q_\text{POD}(t_k;\mathbf p^*) = \mathbf V(\mathbf p^*)\mathbf q(t_k;\mathbf p^{*})$

\ \ \ \ \ \ \ \ \ \ \textbf{End}

Output: $\mathbf Q_\text{POD}(t_k;\mathbf p^{*})$.
\end{itemize}
\caption{POD-PROM framework for linear system~\eqref{eq:4-1HFM}}
\label{alg:pod-prom}
\end{algorithm}

\renewcommand\refname{Reference}
\bibliography{DMD_PROM}

\newpage

\end{document}